\documentclass[article]{amsart}
\usepackage{cases}
\usepackage{amsmath, amsfonts, pifont, amssymb}
\usepackage{verbatim}
\usepackage{geometry}
\newtheorem{theorem}{Theorem}[section]
\newtheorem{lemma}[theorem]{Lemma}

\newtheorem{corollary}[theorem]{Corollary}

\newcommand{\beq}{\begin{equation}}
\newcommand{\eeq}{\end{equation}}
\newcommand{\beqq}{\begin{equation*}}
\newcommand{\eeqq}{\end{equation*}}

\theoremstyle{definition}
\newtheorem{definition}[theorem]{Definition}

\theoremstyle{remark}

\numberwithin{equation}{section}

\numberwithin{equation}{section}

\begin{document}

\title[Scattering for defocusing NLS on $ \mathbb{R}^m$ $\times$ $\mathbb{T}$]{On scattering for the defocusing nonlinear Schr{\"o}dinger equation on waveguide $ \mathbb{R}^m$ $\times$ $\mathbb{T}$ (when $m=2,3$)}

\author{Zehua Zhao}

\maketitle

\begin{abstract}
In the article, we prove the large data scattering for two problems, i.e. the defocusing quintic nonlinear Schr{\"o}dinger equation on $\mathbb{R}^2$ $\times$ $\mathbb{T}$ and the defocusing cubic nonlinear Schr{\"o}dinger equation on $\mathbb{R}^3$ $\times$ $\mathbb{T}$. Both of the two equations are mass supercritical and energy critical. The main ingredients of the proofs contain global Stricharz estimate, profile decomposition and energy induction method. This paper is the second project of our series work (two papers, together with [38]) on large data scattering for the defocusing critical NLS with integer index nonlinearity on low dimensional waveguides. At this point, this category of problems are almost solved except for two remaining resonant system conjectures and the quintic NLS problem on $\mathbb{R}\times \mathbb{T}$. 
\end{abstract}
\bigskip

\bigskip
\section{Introduction}
\noindent First, we consider the following defocusing nonlinear equation with power-type nonlinearity on the waveguide $\mathbb{R}^m \times \mathbb{T}^n$ (the product spaces in this form are called waveguide manifold, see [27, 28] for information) as follows:
\begin{equation}\label{equation}
\aligned
(i\partial_t+ \Delta_{\mathbb{R}^{m} \times \mathbb{T}^{n}}) u &= F(u) = |u|^{p-1} u, \\
u(0,x) &= u_{0} \in H^{1}(\mathbb{R}^{m} \times \mathbb{T}^{n}).
\endaligned
\end{equation}
\noindent where $\Delta_{\mathbb{R}^m\times \mathbb{T}^n}$ is the Laplace-Beltrami operator on $\mathbb{R}^m\times \mathbb{T}^n$ and $u:\mathbb{R}\times \mathbb{R}^m\times \mathbb{T}^n \rightarrow \mathbb{C}$ is a complex-valued function. Specifically, when $m=3,n=1,p=3$, the equation would become cubic NLS on $\mathbb{R}^{3} \times \mathbb{T}$ as follows:
\begin{equation}\label{equation}
\aligned
(i\partial_t+ \Delta_{\mathbb{R}^3\times \mathbb{T}}) u &= F(u) = |u|^{2} u, \\
u(0,x) &= u_{0} \in H^{1}(\mathbb{R}^3\times \mathbb{T}).
\endaligned
\end{equation}
\noindent In addition, when $m=2,n=1,p=5$, the equation would become quintic NLS on $\mathbb{R}^{2} \times \mathbb{T}$ as follows:
\begin{equation}\label{equation}
\aligned
(i\partial_t+ \Delta_{\mathbb{R}^2\times \mathbb{T}}) u &= F(u) = |u|^{4} u, \\
u(0,x) &= u_{0} \in H^{1}(\mathbb{R}^2\times \mathbb{T}).
\endaligned
\end{equation}
\noindent In this paper, we consider the large data scattering for the initial value problem (1.2) and the initial value problem (1.3). There are several reasons why we consider the two problems together and we will explain them shortly. A brief overview of existing related results is as follows:\vspace{3mm}

\noindent For Euclidean case, the large data theory of critical and subcritical NLS is much better understood, at least in the defocusing case (see for example [6, 9, 12, 26, 29]). For waveguide case, we are also naturally interested in the range of the nonlinearity index $p$ in (1.1) when the initial value problem is global well-posed and scattering. According to the existing results and theories, we expect that the solution of (1.1) globally exists and scatters in the range $1+\frac{4}{m}\leq p \leq 1+\frac{4}{m+n-2}$. And for those two problems, the index ($p=3$) in equation (1.2) lies in the range ($\frac{7}{3} \leq p=3\leq 3$); also, the index ($p=5$) in equation (1.3) lies in the range ($3 \leq p=5\leq 5$). So it is reasonable for us to consider those two problems and expect the solutions of (1.2) and (1.3) are global well-posed and scattering.\vspace{3mm}

\noindent There are many related results about NLS on waveguides such as [7] (defocusing cubic $\mathbb{R}^2 \times \mathbb{T}$), [14] (defocusing quintic $\mathbb{R} \times \mathbb{T}^2$), [20] (defocusing cubic $\mathbb{R}\times \mathbb{T}^3$), [36] (defocusing NLS on $\mathbb{R}^d\times \mathbb{T}$) and [38] (defocusing cubic $\mathbb{R}^2\times \mathbb{T}^2$).\vspace{3mm} 

\noindent In [7], X. Cheng, Z. Guo, K. Yang and L. Zhao proved the large data scattering for the the defocusing cubic NLS on $\mathbb{R}^2 \times \mathbb{T}$. One remarkable point is that they have used resonant system approximation. Additionally, in [37], K. Yang and L. Zhao have proved the large data scattering for the corresponding resonant system. \vspace{3mm} 

\noindent In [14], Z. Hani and B. Pausader proved the large data scattering for defocusing quintic NLS on $\mathbb{R} \times \mathbb{T}^2$ based on an assumption, i.e. the large data scattering for a corresponding quintic resonant system. A special case of the quintic resonant system coincides with the 1 dimensional mass-critical NLS problem. The large data scattering for 1 dimensional mass-critical NLS problem is proved by B. Dodson ([11]). Additionally, the way to construct global-in-time Stricharz estimate in this paper is very important and useful.\vspace{3mm}

\noindent In [20], A. Ionescu and B. Pausader proved the global well-posedness of defocusing cubic NLS on $\mathbb{R} \times \mathbb{T}^3$. We are also largely inspired by the method of this paper.\vspace{3mm}

\noindent In [35], N. Tzvetkov and N. Visciglia proved the global well-posedness and scattering for defocusing NLS on $\mathbb{R}^d \times \mathbb{T}$ with nonlinearity index $p$ satisfying $1+\frac{4}{d}< p < 1+\frac{4}{d-1}$. In the case, the equation is mass supercritical and energy subcritical. Specifically, when $d=2$, the corresponding index $p$ satisfies $3<p<5$; when $d=3$, the corresponding index $p$ satisfies $\frac{7}{3}<p<3$. In this paper, we discuss the critical cases when $d=3,p=3$ (equation (1.2)) and $d=2,p=5$ (equation (1.3)). \vspace{3mm}

\noindent In [38], we proved the large data scattering for defocusing cubic NLS on $\mathbb{R}^2 \times \mathbb{T}^2$ based on an assumption, i.e. the large data scattering for a corresponding cubic resonant system. A special case of the quintic resonant system coincides with the 2 dimensional mass-critical NLS problem. The large data scattering for 2 dimensional mass-critical NLS problem is proved by B. Dodson ([10]). Additionally, one dimensional and higher order dimensional mass-critical NLS problems are also solved by B. Dodson ([12]). We also refer to [24, 33, 34] and they are important results about mass-critical NLS.\vspace{3mm}

\noindent Now we consider a series of more specific problems of (1.1), i.e. large data scattering for the defocusing critical NLS with integer index nonlinearity on low dimensional (when $m+n \leq 4$) waveguides. First, noticing the range $1+\frac{4}{m}\leq p \leq 1+\frac{4}{m+n-2}$, we have $n=0,1,2$. When $n=0$, we have $m=4,p=3$ (cubic NLS on $\mathbb{R}^4$), $m=3,p=5$ (quintic NLS on $\mathbb{R}^3$), $m=2,p=3$ (cubic NLS on $\mathbb{R}^2$) , $m=4,p=2$ (4d mass critical NLS) and $m=1,p=5$ (quintic NLS on $\mathbb{R}$). When $n=1$, we have $m=1,p=5$ (quintic NLS on $\mathbb{R}\times \mathbb{T}$),  $m=2,p=5$ (quintic NLS on $\mathbb{R}^2\times \mathbb{T}$, discussed in this paper), $m=3,p=3$ (cubic NLS on $\mathbb{R}^3\times \mathbb{T}$, discussed in this paper) and $m=2, p=3$ (cubic NLS on $\mathbb{R}^2 \times \mathbb{T}$). When $n=2$, we have $m=1,p=5$ (quintic NLS on $\mathbb{R} \times \mathbb{T}^2$) and $m=2,p=3$ (cubic NLS on $\mathbb{R}^2 \times \mathbb{T}^2$). There are totally $11$ specific problems. First, critical NLS problems ($5$ problems) on pure Euclidean domains are well known (see [9, 10, 11, 26]). Additionally, quintic NLS on $\mathbb{R} \times \mathbb{T}$ are expected to be similar to cubic NLS on $\mathbb{R}^2 \times \mathbb{T}$ (see [7]) since both of them are mass critical and energy subcritical. There are $5$ problems left, they are discussed in [7, 14, 38] and this paper.\vspace{3mm}

\noindent Moreover, [2, 6, 15, 16, 17, 18, 21, 22, 25, 26, 29] are some other important resources and related results. Generally speaking, the difficulty of the critical NLS problems on waveguides $\mathbb{R}^m \times \mathbb{T}^n$ increase if the whole dimension $m+n$ increase or if the $\mathbb{R}$-dimension $m$ decrease. As for the introduction of the related NLS problems on waveguides, we also refer to [7, Introduction], [14, Introduction], [20, Introduction] and [38, Introduction] for more information.\vspace{3mm}

\noindent \textbf{A word on quintic 3d problems:} While scattering holds for the quintic equation on $\mathbb{R}^3$ (see [9]), $\mathbb{R}^2\times \mathbb{T}$ (see this paper) and $\mathbb{R}\times \mathbb{T}^2$ (see [14]), it is not expected to hold on $\mathbb{T}^3$. The situation on $\mathbb{R}\times \mathbb{T}^2$ seems to be a borderline case for this question, i.e. defocusing quintic NLS equation on three dimensional waveguides. \vspace{3mm}

\noindent \textbf{A word on cubic 4d problems:} Also, while scattering holds for the cubic equation on $\mathbb{R}^4$ (see [26]), $\mathbb{R}^3\times \mathbb{T}$ (see this paper) and $\mathbb{R}^2\times \mathbb{T}^2$, it is not expected to hold on $\mathbb{R}\times \mathbb{T}^3$ or $\mathbb{T}^4$. The situation on $\mathbb{R}^2\times \mathbb{T}^2$ seems to be a borderline case for this question, i.e. defocusing cubic NLS equation on four dimensional waveguides. \vspace{3mm}

\noindent \textbf{A word on the remaining problems in this category: }Together with [7, 14, 38] (corresponding to cubic $\mathbb{R}^2\times \mathbb{T}$ problem, quintic $\mathbb{R}\times \mathbb{T}^2$ problem and cubic $\mathbb{R}^2\times \mathbb{T}^2$ problem respectively), \emph{large data scattering for defocusing critical NLS with integer-index nonlinearity on low dimensional waveguides} are almost solved except for the two resonant system conjectures arising from [14, 38] and the quintic NLS on $\mathbb{R}\times \mathbb{T}$. One inspiring thing is, in [37], the corresponding cubic resonant system conjecture (arises from cubic $\mathbb{R}^2\times \mathbb{T}$ problem) is solved by K. Yang and L. Zhao. We expect one may use [37] and the large data scattering results for defocusing mass critical NLS by B. Dodson ([10, 11, 12]) to prove the two remaining resonant system conjectures. Additionally, we expect one may learn from [7] and use the scattering result for 1d mass critical NLS problem by B. Dodson ([11]) to prove the large data scattering for the quintic NLS on $\mathbb{R}\times \mathbb{T}$. In some sense, quintic NLS problem on $\mathbb{R}\times \mathbb{T}$ is the quintic analogue of cubic NLS problem on $\mathbb{R}^2\times \mathbb{T}$ ([7]).\vspace{3mm}

\noindent \textbf{Comparison of the two problems:} There are some similarities and differences between the two problems. One significant \emph{similarity} is that both of the two equations are mass supercritical and energy critical which leads that the spirits of the linear profiles and the profile decompositions are same. As for the \emph{differences}, (1.2) is cubic and (1.3) is quintic thus the nonlinear estimates would be different. Another difference of those two problems would be the whole spatial dimensions ($3$ and $4$). Those two equations are analogues of each other in some sense.\vspace{3mm}

\noindent The following two Theorems are the main results of this paper. Theorem 1.1 is for equation (1.2) (cubic $\mathbb{R}^3\times \mathbb{T}$ problem) and Theorem 1.2 is for equation (1.3) (quintic $\mathbb{R}^2\times \mathbb{T}$ problem).

\begin{theorem}For any initial data $u_0\in H^1(\mathbb{R}^3\times \mathbb{T})$, there exists a solution $u\in X^1_{c}(\mathbb{R})$ to (1.2) that is global and scattering in the sense that there exists $v^{\pm \infty} \in H^1(\mathbb{R}^3\times \mathbb{T})$ such that
\begin{equation}
||u(t)-e^{it\Delta_{\mathbb{R}^3\times \mathbb{T}}}v^{\pm \infty}||_{H^1(\mathbb{R}^3\times \mathbb{T})} \rightarrow 0,\textmd{ as } t\rightarrow \pm \infty.
\end{equation}
\end{theorem}
\noindent The uniqueness space $X^1_c\subset C_t(\mathbb{R}:H^1(\mathbb{R}^3\times \mathbb{T}))$ was first introduced by Herr-Tataru-Tzvetkov ([16]) (see also [14], we define the solution space in a similar way.). 

\begin{theorem}For any initial data $u_0\in H^1(\mathbb{R}^2\times \mathbb{T})$, there exists a solution $u\in X^1_{c}(\mathbb{R})$ to (1.3) that is global and scattering in the sense that there exists $v^{\pm \infty} \in H^1(\mathbb{R}^2\times \mathbb{T})$ such that
\begin{equation}
||u(t)-e^{it\Delta_{\mathbb{R}^2\times \mathbb{T}}}v^{\pm \infty}||_{H^1(\mathbb{R}^2\times \mathbb{T})} \rightarrow 0,\textmd{ as } t\rightarrow \pm \infty
\end{equation}
\end{theorem}
\noindent The uniqueness space $X^1_c\subset C_t(\mathbb{R}:H^1(\mathbb{R}^2\times \mathbb{T}))$ was first introduced by Herr-Tataru-Tzvetkov ([16]) (see also [14], we define the solution space in a similar way.).\vspace{3mm} 

\noindent For those two problems, we will prove small data result first by using global Stricharz estimate and some standard arguments and then we can use profile decomposition and energy induction method to extend our analysis to large data case.\vspace{3mm}

\noindent Similar as other related results (see [7, 14, 20, 38]), profile decomposition is a crucial step. In order to understand the appearance of the profiles, specifically for this problem, in view of the scaling-invariant of (1.1) under 
\[\mathbb{R}_x^3 \times \mathbb{T}_y \rightarrow M_{\lambda}:=\mathbb{R}_x^3 \times (\lambda^{-1}\mathbb{T})_{y},\quad u \rightarrow \tilde{u}(x,y,t)=\lambda u(\lambda x,\lambda y,\lambda^2 t).\]

\noindent \emph{Remark.} Similar scaling-invariance analysis also works for the case of quintic NLS on $\mathbb{R}^2 \times \mathbb{T}$ so we omit it.\vspace{3mm}

\noindent There are two extreme situations as follows:\vspace{3mm}

\noindent When $\lambda \rightarrow 0$, the manifolds $M_\lambda$ will be similar to $\mathbb{R}^4$ and we can use the scattering result for four dimensional energy critical NLS by E. Ryckman and M. Visan ([26]). The appearance is a manifestation of the energy-critical nature of the nonlinearity. This corresponds in $M_1$ to solutions with initial data (This behavior corresponds to Euclidean profiles which we will analyze more precisely in section 5 and section 6.)
\[ u^{\lambda}(x,y,0)=\lambda^{-1}\phi(\lambda^{-1} (x,y)) \quad \phi \in C^{\infty}_0(\mathbb{R}^4),\lambda \rightarrow 0.
\]

\noindent When $\lambda \rightarrow \infty$, the manifolds $M_\lambda$ become thinner and thinner and resemble $\mathbb{R}^3$. The problem will become similar to the following cubic NLS problem on $\mathbb{R}^3$ :
\[
(i\partial_t+\Delta_x )u=|u|^2u ,\quad u(0) \in H^1(\mathbb{R}^3) .
\]
\noindent Those solutions on $M_{\lambda}$ correspond to solutions on $M_1$ with initial data
\begin{equation}
u^{\lambda}(x,y,0)=\lambda^{-1}\phi(\lambda^{-1}x,y) \quad \phi \in C^{\infty}_0(\mathbb{R}^3\times \mathbb{T}),\lambda \rightarrow \infty.
\end{equation}
\noindent Different from [14] (defocusing quintic $\mathbb{R}\times \mathbb{T}^2$ problem) and [38] (defocusing cubic $\mathbb{R}^2\times \mathbb{T}^2$ problem), for those two problems in this paper, large scale profile will not appear in the profile decomposition. Put in another way, extracting orthogonal Euclidean profiles and scale-one profiles are enough for us to control the scattering norm of the linear Schr{\"o}dinger propagation of the remainder flow and extend our analysis to large data case. The reason is: in (1.6) the $\dot{H}^1$ norm of $u^{\lambda}$ will converge to $0$ when $\lambda \rightarrow \infty$; additionally, the $H^1$ boundedness condition ensures the mass of $\phi$ is $0$, otherwise the mass of $u^{\lambda}$ will blow up. Thus $u^{\lambda}$ must converge to $0$ in $H^1$. \vspace{3mm}

\noindent Both of the domain and the nonlinearity in the equation play important roles in the asymptotic behavior of the solutions when we consider NLS problems. That is why there are some remarkable differences among those related problems. In [14, 38] (the equation is both mass critical and energy subcritical), there are three types of profiles considered, i.e. Euclidean profiles, scale-one profiles and large-scale profiles. In [20] and this paper (the equations are only energy critical), Euclidean profiles and scale-one profiles are considered. In [7] (the equation is mass critical and energy subcritical), large-scale profiles and scale-one profiles are considered.\vspace{3mm}

\noindent The proofs of the Theorem 1.1 and Theorem 1.2 follow from a standard skeleton based on the Kenig-Merle machinery [21, 22]. Mainly there are three important ingredients: global Strichartz estimates, profile decomposition and energy induction method.\vspace{3mm}

\noindent \textbf{The organization of this paper}: in Section 2, we introduce some notations and function spaces; in Section 3, we prove the global Strichartz estimates; in Section 4, we prove the local well-posedness and small data scattering of (1.2) and (1.3); in Section 5, we describe Euclidean profile and Euclidean approximation; in Section 6, we obtain the linear profile decomposition that leads us to analyze the large data case; in Section 7, we prove the contradiction argument leading to our large data scattering result, i.e. Theorem 1.1 and Theorem 1.2. Eventually, \emph{an important claim is:} since the two problems have many similarities, for convenience, when possible, we will mainly discuss cubic $\mathbb{R}^3\times \mathbb{T}$ problem and when necessary we will clarify the differences with details.\vspace{3mm}  

\section{Notation and Preliminaries}
\noindent About the notation, we write $A \lesssim B$ to say that there is a constant $C$ such that $A\leq CB$. We use $A \simeq B$ when $A \lesssim B \lesssim A $. Particularly, we write $A \lesssim_u B$ to express that $A\leq C(u)B$ for some constant $C(u)$ depending on $u$.\vspace{3mm}

\noindent We also define the partial Littlewood-Paley projectors $P^x_{\leq N}$ and $P^x_{\geq N}$ as follows: fix a real-valued radially symmetrically bump function $\varphi(\xi)$ satisfying
\begin{align*}
\varphi(\xi)=
\begin{cases}
1, |\xi|\leq 1, \\
0, |\xi| \geq 2,
\end{cases}
\end{align*}
\noindent for any dyadic number $N\in 2^{\mathbb{Z}}$, let
\[ \mathcal{F}_x(P^x_{\leq N} f)(\xi,y)=\varphi(\frac{\xi}{N})(\mathcal{F}_xf)(\xi,y),
\]
\[\mathcal{F}_x(P^x_{\geq N} f)(\xi,y)=(1-\varphi(\frac{\xi}{N}))(\mathcal{F}_xf)(\xi,y).
\]
\noindent \textbf{Function spaces.} In this paper, we use some function spaces ($X^s,Y^s,N^s$) based on atomic space and variation space. Those spaces were essentially introduced in [16]. (see also [17]). Moreover, we construct those spaces as in [14, 38]. We also use $X^1_c$ to be the basic solution space. Those spaces have some nice properties. For convenience, we give some basic definitions for the case ``cubic $\mathbb{R}^3\times \mathbb{T}$'' as follows and and we refer to [16, 17] for more description of those spaces. \vspace{3mm}

\noindent \emph{Remark. }The setting for quintic $\mathbb{R}^2\times \mathbb{T}$ problem is similar and we omit it. \vspace{3mm}

\noindent For $C=[-\frac{1}{2},\frac{1}{2})^4 \in \mathbb{R}^4$ and $z\in \mathbb{R}^4$, we denote by $C_z=z+C$ the translate by $z$ and define the projection operator $P_{C_z}$ as follows, ($\mathcal{F}$ is the Fourier transform): 
\[
\mathcal{F}(P_{C_z} f)=\chi_{C_z}(\xi) \mathcal{F} (f)  (\xi).
\]\noindent As in [16, 17], for $s\in \mathbb{R}$, we define:
\[ \|u\|_{X_0^s(\mathbb{R})}^2=\sum_{z\in \mathbb{Z}^4} \langle z \rangle^{2s} \|P_{C_z} u\|_{U_{\Delta}^2(\mathbb{R};L^2)}^2
\]
\noindent and similarly we have,  
\[\|u\|_{Y^s(\mathbb{R})}^2=\sum_{z\in \mathbb{Z}^4} \langle z \rangle^{2s} \|P_{C_z} u\|_{V_{\Delta}^2(\mathbb{R};L^2)}^2\]          
\noindent where the $U_{\Delta}^p$ and $V_{\Delta}^p$ are the atomic and variation spaces respectively of functions on $\mathbb{R}$ taking values in $L^2(\mathbb{R}^3\times \mathbb{T})$.\vspace{3mm}

\noindent Moreover, for an interval $I \subset \mathbb{R}$, we can also define the restriction norms $X^s_0(I)$ and $Y^s(I)$ in the natural way:
\noindent                                               $||u||_{X_0^s(I)}=$ inf $\{||v||_{X_0^s(\mathbb{R})}:v\in X_0^s(\mathbb{R})$ satisfying $v_{|I}=u_{|I}\}$. And similarly for $Y^s(I)$. Additionally, a modification for to $X_0^s(\mathbb{R})$: 

\noindent $X^s(\mathbb{R}):=\{u:\phi_{-\infty}=\lim\limits_{t\rightarrow -\infty}e^{-it\Delta}u(t) \textmd{ exists in }  H^s,u(t)-e^{it\Delta}\phi_{-\infty} \in X^s_0(\mathbb{R})\} $ equipped with the norm:
\begin{equation}
||u||^2_{X^s(\mathbb{R})}=||\phi_{-\infty}||^2_{H^s(\mathbb{R}^3\times \mathbb{T})}+||u-e^{it\Delta}\phi_{-\infty}||^2_{X_0^s(\mathbb{R})}.
\end{equation}\noindent Our basic space to control solutions is $X^1_c(I)=X^1(I)\cap C(I:H^1).$  Also we use $X^1_{c,loc}(I)$ to express the set of all solutions in $C_{loc}(I:H^1)$ whose $X^1(J)$-norm is finite for any compact subset $J \subset I$.\vspace{3mm}

\noindent At last, in order to control the nonlinearity on interval $I$, we need to define `$N$ -Norm' as follows, on an interval $I=(a,b)$ we have:
\begin{equation}\label{equation}
\| h\|_{N^s(I)}=\|\int_{a}^{t} e^{i(t-s)\Delta} h(s) ds \|_{X^s(I)} .
\end{equation}
\noindent We also need the following theorem which has analogues in [14, 16, 17].
\begin{theorem}\label{theorem} If $f\in L^1_t(I,H^1(\mathbb{R}^3\times \mathbb{T}))$, then 
\[||f||_{N(I)} \lesssim \sup_{v\in Y^{-1}(I),||v||_{Y^{-1}(I)}\leq 1} \int_{I \times (\mathbb{R}^3 \times \mathbb{T})} f(x,t)\overline{v(x,t)}dxdt.
\]
Also, we have the following estimate holds for any smooth function $g$ on an interval $I=[a,b]$: 
\[ ||g||_{X^1(I)}\lesssim ||g(0)||_{H^1(\mathbb{R}^3\times \mathbb{T})}+(\sum_N ||P_N(i\partial_t+\Delta)g||^2_{L^1_t(I,H^1(\mathbb{R}^3\times \mathbb{T}))})^{\frac{1}{2}}.
\]
\end{theorem}
\noindent \emph{Remark.} We have analogue of the above lemma for quintic $\mathbb{R}^2\times \mathbb{T}$ problem.\vspace{3mm}

\noindent The following lemma is also useful. It is exactly the analogue of [20, Lemma 2.3].
\begin{lemma}
\noindent For $f \in H^1(\mathbb{R}^3 \times \mathbb{T})$, there holds that
\begin{equation}
||f||_{L^4(\mathbb{R}^3 \times \mathbb{T})}\lesssim (\sup_N N^{-1}||P_N f||_{L^{\infty}(\mathbb{R}^3 \times \mathbb{T})})^{\frac{1}{2}}(||f||_{H^1(\mathbb{R}^3 \times \mathbb{T})})^{\frac{1}{2}}.
\end{equation}
\end{lemma}
\section{Global Strichartz estimate}
\noindent \textbf{For cubic $\mathbb{R}^3\times \mathbb{T}$ problem:} we have
\begin{theorem}\label{theorem}Now we prove the following Strichartz Estimate:
\begin{equation}\label{equation}
\|e^{it\Delta_{\mathbb{R}^3\times\mathbb{T}}}P_{\leq N}u_0\|_{{l_{\gamma}^q L_{x,y,t}^p}(\mathbb{R}^3\times\mathbb{T}\times[2\pi\gamma,2\pi(\gamma+1)])}\lesssim N^{2-\frac{6}{p}}\|u_0\|_{L^2( \mathbb{R}^3\times\mathbb{T})} 
\end{equation}
\noindent whenever 
\begin{equation}\label{equation} 
\frac{22}{7}<p<6 \quad and \quad \frac{2}{q}+\frac{3}{p}=\frac{3}{2}.
\end{equation}
\end{theorem}
\noindent \emph{Proof:} The main idea of the proof is similar to [14, Theorem 3.1] (see also [36, Theorem 3.1]). We use duality argument, $T-T^*$ argument, a partition of unity and then estimate the diagonal part and non-diagonal part separately to obtain the above estimate. First, let us prove a more precise conclusion and we can get the above estimate by duality:

\begin{lemma}
For any $h\in C_c^{\infty}(\mathbb{R}_x^3\times \mathbb{T}_y \times \mathbb{R}_t)$, there holds that 
\begin{equation}\label{equation}
\aligned
&\|\int_{s\in \mathbb{R}} e^{-is\Delta_{\mathbb{R}^3\times \mathbb{T}}} P_{\leq N} h(x,y,s) ds\|_{L_{x,y}^{2}(\mathbb{R}_x^3\times \mathbb{T}_y)} \\ &\lesssim N^{2-\frac{6}{p}} \|h\|_{l_{\gamma}^2 L_{x,y,t}^{p^{'}}(\mathbb{R}^3 \times \mathbb{T} \times [2\pi\gamma,2\pi(\gamma+1)])}+N^{1-\frac{20}{7p}} \|h\|_{l_{\gamma}^{q^{'}} L_{x,y,t}^{p^{'}}(\mathbb{R}^3\times \mathbb{T} \times [2\pi\gamma,2\pi(\gamma+1)])}  
\endaligned
\end{equation}
for any (p,q) satisfies (3.2).
\end{lemma}
\noindent \emph{Remark.} The above threshold $p<6$ results from the requirement $q^{'}< 2$ noticing the inclusion property of the sequence spaces.
\begin{corollary}
Noticing the sequence space inclusion and Lemma 3.2, the following Strichartz estimate also holds:
\begin{equation}\label{equation}
\|e^{it\Delta_{\mathbb{R}^3\times\mathbb{T}}}P_{\leq N}u_0\|_{{l_{\gamma}^q L_{x,y,t}^p}(\mathbb{R}^3\times\mathbb{T}\times[2\pi\gamma,2\pi(\gamma+1)])}\lesssim N^{2-\frac{6}{p}}\|u_0\|_{L^2( \mathbb{R}^3\times\mathbb{T})} 
\end{equation}
\noindent whenever 
\begin{equation}\label{equation} 
p>\frac{22}{7} \quad and \quad \frac{2}{q}+\frac{3}{p}=1.
\end{equation}
\end{corollary}
\noindent \emph{Remark.} This estimate corresponds to the scattering norm which we will define in Section 4.\vspace{3mm}

\noindent \emph{Proof of Lemma 3.2:} In order to distinguish between the large and small time scales, we choose a smooth partition of unity $1=\sum\limits_{\gamma\in \mathbb{Z}} \chi(t-2\pi\gamma)$ with $\chi$ supported in $[-2\pi,2\pi]$. We also denote by $h_{\alpha}(t)=\chi(t) h(2\pi\alpha+t)$. Using the semigroup property and the unitarity of $e^{it\Delta_{\mathbb{R}^3\times \mathbb{T}}}$ we can get :

\begin{align*}
&\|\int_{s\in \mathbb{R}} e^{-is\Delta_{\mathbb{R}^3\times \mathbb{T}}} P_{\leq N} h(x,y,s)ds\|_{L_{x,y}^2(\mathbb{R}^3\times \mathbb{T})}\\
&=\int_{s,t\in \mathbb{R}} \langle e^{-is\Delta_{\mathbb{R}^3\times \mathbb{T}}} P_{\leq N} h(s),e^{-it\Delta_{\mathbb{R}^3\times \mathbb{T}}} P_{\leq N} h(t) \rangle_{L_{x,y}^2(\mathbb{R}^3\times \mathbb{T})\times L_{x,y}^2(\mathbb{R}^3\times \mathbb{T})}dsdt  \\
&=\sum_{\alpha,\beta} \int_{s,t\in \mathbb{R}} \langle \chi(s-2\pi\alpha) e^{-is\Delta_{\mathbb{R}^3\times \mathbb{T}}} P_{\leq N} h(s),\chi(s-2\pi\beta) e^{-is\Delta_{\mathbb{R}^3\times \mathbb{T}}} P_{\leq N} h(t)\rangle_{L_{x,y}^2(\mathbb{R}^3\times \mathbb{T})\times L_{x,y}^2(\mathbb{R}^3\times \mathbb{T})}dsdt \\
&=\sum_{\alpha,\beta} \int_{s,t\in [-2\pi,2\pi]}\langle e^{-i(2\pi(\alpha-\beta)+s))\Delta_{\mathbb{R}^3\times \mathbb{T}}} P_{\leq N}h_{\alpha}(s),e^{-it\Delta_{\mathbb{R}^3\times \mathbb{T}}} P_{\leq N}h_{\beta}(t)\rangle_{L_{x,y}^2 \times L_{x,y}^2} dsdt \\
&=\sigma_{d}+\sigma_{nd}.
\end{align*}
\noindent Here we have,
\[
\sigma_{d}=\sum_{\alpha\in \mathbb{Z},| \gamma| \leq 9} \int_{s,t\in \mathbb{R}} \langle e^{-i(s-2\pi\gamma)\Delta_{\mathbb{R}^3\times \mathbb{T}}} P_{\leq N}h_{\alpha}(s),e^{-it\Delta_{\mathbb{R}^3\times \mathbb{T}}} P_{\leq N}h_{\alpha+\gamma}(t) \rangle_{L_{x,y}^2 \times L_{x,y}^2} dsdt.
\]
\[
\sigma_{nd}=\sum_{\alpha,\gamma \in \mathbb{Z},| \gamma| > 10} \int_{s,t\in \mathbb{R}} \langle e^{-i(s-2\pi\gamma)\Delta_{\mathbb{R}^3\times \mathbb{T}}} P_{\leq N}h_{\alpha}(s),e^{-it\Delta_{\mathbb{R}^3\times \mathbb{T}}} P_{\leq N}h_{\alpha+\gamma}(t) \rangle _{L_{x,y}^2 \times L_{x,y}^2} dsdt.
\]
\noindent Here, `$d$' is short for `diagonal' and `$nd$' is short for `non-diagonal'. We will estimate the diagonal part and the non-diagonal part by using different methods.\vspace{3mm}

\noindent \emph{For the diagonal part:} First, we need to use a local-in-time $L^p$ estimate as follows:
\begin{lemma} 
\noindent Let $p_1=\frac{22}{7}$, then for any $p>p_1$, $N\geq 1$,and $f\in L^2(\mathbb{R}^3 \times \mathbb{T})$,
\begin{equation}
||e^{it\Delta}P_N f||_{L^p(\mathbb{R}^3 \times \mathbb{T} \times [0,2\pi])}\lesssim_p N^{2-\frac{6}{p}}||f||_{L^2(\mathbb{R}^3\times \mathbb{T})}.
\end{equation}
\end{lemma}

\noindent The proof of the local-in-time estimate can also be proved in a similar way as in [20] with small modifications (see also [25] which is a related recent result). We omit it.\vspace{3mm}

\noindent According to the estimate (3.6) above, by duality we have
\begin{equation}
\|\int_{s\in \mathbb{R}} e^{-is\Delta_{\mathbb{R}^2\times \mathbb{T}^2}} P_{\leq N} h(s) ds\|_{L_{x,y}^2(\mathbb{R}^3\times \mathbb{T})} \lesssim N^{2-\frac{6}{p}} \|h\|_{L_{x,y,t}^{p^{'}}(\mathbb{R}^3\times \mathbb{T}\times [-2\pi,2\pi])} 
\end{equation}
\noindent where $h$ is supported in $[-2\pi,2\pi]$. And consequently,
\begin{equation}\label{equation}
\aligned
\sigma_{d} &=\sum_{\alpha \in \mathbb{Z} ,|\gamma| \leq 9} \int_{s,t\in \mathbb{R}} \langle e^{-i(s-2\pi\gamma)\Delta_{\mathbb{R}^3\times \mathbb{T}}} P_{\leq N}h_{\alpha}(s),e^{-it\Delta_{\mathbb{R}^3\times \mathbb{T}}} P_{\leq N}h_{\alpha+\gamma}(t)\rangle_{L_{x,y}^2 \times L_{x,y}^2} dsdt  \\
&\leq \sum_{\alpha\in \mathbb{Z} ,|\gamma| \leq 9} \|\int_{s\in \mathbb{R}} e^{-is\Delta_{\mathbb{R}^3 \times \mathbb{T}}} P_{\leq N} h_{\alpha}(2\pi\gamma+s)ds\|_{L_{x,y}^2} \|\int_{s\in \mathbb{R}} e^{-is\Delta_{\mathbb{R}^3 \times \mathbb{T}}} P_{\leq N} h_{\alpha+\gamma}(s)ds\|_{L_{x,y}^2} \\
&\lesssim N^{2(2-\frac{6}{p})} \sum_{\alpha} \|h_{\alpha}\|_{L_{x,y,t}^{p^{'}}(\mathbb{R}^3 \times \mathbb{T} \times [-2\pi,2\pi])}^2.  
\endaligned
\end{equation}

\noindent This finishes the estimate for the diagonal part in (3.3).\vspace{3mm}

\noindent \emph{For the non-diagonal part:} we need to use Lemma 3.5 that we will prove soon and we can apply it to estimate the non-diagonal part by using H$\ddot{o}$lder's inequality and the discrete Hardy-Sobolev inequality as below:
\begin{equation}\label{equation}
\aligned
\sigma_{nd}&=\sum_{\alpha,\gamma \in \mathbb{Z},| \gamma| > 10} \int_{t \in \mathbb{R}} \langle \int_{s \in \mathbb{R}}e^{-i(s-2\pi\gamma)\Delta_{\mathbb{R}^3\times \mathbb{T}}} P_{\leq N}h_{\alpha}(s)ds,e^{-it\Delta_{\mathbb{R}^3\times \mathbb{T}}} P_{\leq N}h_{\alpha+\gamma}(t) \rangle_{L_{x,y}^2 \times L_{x,y}^2} dt \\
&\lesssim N^{2-\frac{40}{7p}} \sum_{\alpha,\gamma \in \mathbb{Z},| \gamma| > 3} |\gamma|^{\frac{3}{p}-\frac{3}{2}} \|h_{\alpha}\|_{L_{x,y,t}^{p^{'}}}\|h_{\alpha+\gamma}\|_{L_{x,y,t}^{p^{'}}}    \\
&\lesssim N^{2-\frac{40}{7p}}\|h_{\alpha}\|_{l_{\alpha}^{q^{'}} L_{x,y,t}^{p^{'}} (\mathbb{R}^3\times \mathbb{T} \times[-2\pi,2\pi])}^2 .
\endaligned
\end{equation}
\begin{lemma} Suppose $\gamma\in \mathbb{Z}$ satisfies $|\gamma|\geq 3$ and that $p> \frac{22}{7}$. For any function $h\in L_{x,y,s}^{p^{'}}(\mathbb{R}^3 \times \mathbb{T} \times [-2\pi,2\pi])$, there holds that:
\[
\|\int_{s\in \mathbb{R}} \chi(s) e^{i(t-s+2\pi\gamma)\Delta_{\mathbb{R}^3 \times \mathbb{T}}} P_{\leq N} h(s)ds\|_{L_{x,y,t}^{p}(\mathbb{R}^3\times \mathbb{T} \times [-2\pi,2\pi])} \lesssim |\gamma|^{\frac{3}{p}-\frac{3}{2}} N^{2-\frac{40}{7p}} \|h\|_{L_{x,y,s}^{p^{'}}(\mathbb{R}^3\times \mathbb{T} \times [-2\pi,2\pi])}.
\]\end{lemma}

\noindent \emph{Proof:} The proof of this lemma is similar to [14, Lemma 3.3] (see also [38]) by using Hardy-Littlewood circle method. The main idea of the proof is to study the Kernel $K_{N,\gamma}$, use a partition and decompose the corresponding index set into three parts and estimate over the three parts separately. \vspace{3mm}

\noindent Without loss of generality, we assume that:
\[h=\chi(s)P_{\leq N}h,\quad ||h||_{L^{p^{'}}(\mathbb{R}^3\times \mathbb{T}\times [-2\pi,2\pi])}=1
\]
\noindent and we define:
\[g(x,y,s)=\int_{s\in \mathbb{R}}e^{i(t-s+2\pi \gamma)\Delta_{\mathbb{R}^3\times \mathbb{T}}} h(x,y,s) ds.
\]
\noindent Also the Kernel is defined as follows: 
\begin{equation}\label{equation}
\aligned
K_N(x,y,t)&=\sum_{k \in \mathbb{Z}}\int_{\mathbb{R}^3_{\xi}}[\eta^1_{\leq N}(\xi_1)\eta^1_{\leq N}(\xi_2)\eta^1_{\leq N}(\xi_3)]^2[\eta^1_{\leq N}(k) ]^2e^{i[x\cdot \xi+y\cdot k+t(|k|^2+|\xi|^2)]} d\xi\\
&=[\int_{\mathbb{R}^3_{\xi}}[\eta^1_{\leq N}(\xi_1) \eta^1_{\leq N}(\xi_2)\eta^1_{\leq N}(\xi_3)]^2e^{i[x\cdot \xi+t|\xi|^2]}]\cdot[\sum_{k \in \mathbb{Z}}[\eta^1_{\leq N}(k) ]^2e^{i[y\cdot k+t|k|^2]}]\\
&=K^{\mathbb{R}^3}_N(x,t)\bigotimes K^{\mathbb{T}}_N(x,t).
\endaligned
\end{equation}
\noindent And we define $K_{N,\gamma}(x,y,t):=K_N(x,y,2\pi\gamma+t)$, and so we have $g(x,y,t)=K_{N,\gamma}*h$. Notice that a remarkable difference is the non-stationary phase estimate because of the dimension, we will have:
\[
\|K_{N,\gamma}\|_{L_{x,y,t}^{\infty}} \lesssim |\gamma|^{-\frac{3}{2}} N
\]
\noindent instead of
\[
\|K_{N,\gamma}\|_{L_{x,y,t}^{\infty}} \lesssim |\gamma|^{-\frac{1}{2}} N^2.
\]
\noindent And
\[||\mathcal{F}_{x,y,t} K_{N,\gamma}||_{L_{\xi,k,\tau}^{\infty}} \lesssim 1
\]
still holds.\vspace{3mm}

\noindent For $\alpha$ a dyadic number, we define $g^{\alpha}(x,y,t)=\alpha^{-1}g(x,y,t)1_{\{\frac{\alpha}{2} \leq |g|\leq \alpha\}}$ which has modulus in $[\frac{1}{2},1]$. We define $h^{\beta}$ similarly for $\beta\in 2^{\mathbb{Z}}$. And we have the following decomposition:
\begin{equation}\label{equation}
\aligned
||g||^p_{L^p_{x,y,t}}&=\langle |g|^{p-2}g,g \rangle \\
&=\langle \sum_{\alpha,\beta}\alpha^{p-1}|\alpha|^{p-2}\alpha,K_{N,\gamma}*h^{\beta} \rangle \\
&=\sum_{\alpha,\beta}\alpha^{p-1} \beta \langle |g^{\alpha}|^{p-2}g^{\alpha},K_{N,\gamma}*h^{\beta} \rangle \\
&=[\sum_{\mathcal{S}_1}+\sum_{\mathcal{S}_2}+\sum_{\mathcal{S}_3}]\alpha^{p-1} \beta \langle |g^{\alpha}|^{p-2}g^{\alpha},K_{N,\gamma}*h^{\beta} \rangle\\
&=\sum_{1}+\sum_{2}+\sum_{3},
\endaligned
\end{equation}
\noindent where $\mathcal{S}_1$, $\mathcal{S}_2$, $\mathcal{S}_3$ are three index sets. Furthermore, we have the following decomposition:
\begin{enumerate}
\item $\mathcal{S}_1=\{(\alpha,\beta):C|\gamma|^{-\frac{3}{2}}N\leq \alpha \beta^{{p^{'}}-1}\}$, 
\item $\mathcal{S}_2=\{(\alpha,\beta):\alpha \beta^{{p^{'}}-1}\leq CN^{\frac{1}{2}}|\gamma|^{-\frac{3}{2}}\}$,
\item $\mathcal{S}_3=\{(\alpha,\beta):CN^{\frac{1}{2}}|\gamma|^{-\frac{3}{2}} \leq \alpha \beta^{{p^{'}}-1} \leq C|\gamma|^{-\frac{3}{2}}N\}$
\end{enumerate}
\noindent for $C$ a large constant to be decided later. For fixed $\alpha,\beta$, we will decompose $K_{N,\gamma}=K^1_{N,\gamma;\alpha,\beta}+K^2_{N,\gamma;\alpha,\beta}$ and estimate them as follows.
\begin{equation}
\langle |g^{\alpha}|^{p-2}g^{\alpha},K^{1}_{N,\gamma}*h^{\beta} \rangle \lesssim ||K^1_{N,\gamma;\alpha,\beta}||_{L^{\infty}_{x,y,t}}||g^{\alpha}||_{L^1}||h^{\beta}||_{L^1},
\end{equation}
\begin{equation}
\langle |g^{\alpha}|^{p-2}g^{\alpha},K^{2}_{N,\gamma}*h^{\beta} \rangle \lesssim ||\mathcal{F}_{x,y,t} K^2_{N,\gamma;\alpha,\beta}||_{L^{\infty}_{\xi,k,\tau}}||g^{\alpha}||_{L^2}||h^{\beta}||_{L^2}.
\end{equation}
\noindent Then we can estimate the three parts as in [14]. The small level and large level are easy to handle, while the estimate for the medium level requires a more delicate decomposition for the kernel as follows:
\begin{equation}
||K^1_{N,\gamma;\alpha,\beta}||_{L^{\infty}_{x,y,t}}\lesssim \alpha \beta^{{p^{'}}-1},
\end{equation}
\begin{equation}
||\mathcal{F}_{x,y,t} K^2_{N,\gamma;\alpha,\beta}||_{L^{\infty}_{\xi,k,\tau}} \lesssim N^{\epsilon} (\alpha \beta^{{p^{'}}-1})^{-1} |\gamma|^{-\frac{3}{2}}.
\end{equation}

\noindent The above kernel decomposition is similar as in [14] by using Hardy-Littlewood circle method. The rest follows as in [14] so we omit it. Eventually the conclusion estimate is as follows:
\[||g||^p_{L^p_{x,y,t}}\lesssim_C ||g||^{\frac{p}{2}}_{L^p_{x,y,t}}|\gamma|^{\frac{6-3p}{4}}\textmd{max}(N^{p-4+\epsilon},N^{\frac{p-2}{4}})\lesssim ||g||^{\frac{p}{2}}_{L^p_{x,y,t}}|\gamma|^{\frac{6-3p}{4}}N^{p-\frac{20}{7}}
\]
\noindent whenever $p>\frac{22}{7}$. That finishes the proof of the non-diagonal estimate.\vspace{3mm}

\noindent \textbf{For quintic $\mathbb{R}^2\times \mathbb{T}$ problem:} we have the following estimate:
\begin{theorem}\label{theorem}Now we prove the following Strichartz Estimate:
\begin{equation}\label{equation}
\|e^{it\Delta_{\mathbb{R}^2\times\mathbb{T}}}P_{\leq N}u_0\|_{{l_{\gamma}^q L_{x,y,t}^p}(\mathbb{R}^2\times\mathbb{T}\times[2\pi\gamma,2\pi(\gamma+1)])}\lesssim N^{\frac{3}{2}-\frac{5}{p}}\|u_0\|_{L^2( \mathbb{R}^2\times\mathbb{T})} 
\end{equation}
\noindent whenever 
\begin{equation}\label{equation} 
p>\frac{18}{5} \quad and \quad \frac{1}{q}+\frac{1}{p}=\frac{1}{2}.
\end{equation}
\end{theorem}
\begin{corollary}According to the inclusion property of sequence spaces, it is not hard to verity the following estimate also holds:
\begin{equation}\label{equation}
\|e^{it\Delta_{\mathbb{R}^2\times\mathbb{T}}}P_{\leq N}u_0\|_{{l_{\gamma}^q L_{x,y,t}^p}(\mathbb{R}^2\times\mathbb{T}\times[2\pi\gamma,2\pi(\gamma+1)])}\lesssim N^{\frac{3}{2}-\frac{5}{p}}\|u_0\|_{L^2( \mathbb{R}^2\times\mathbb{T})} 
\end{equation}
\noindent whenever 
\begin{equation}\label{equation} 
p>\frac{18}{5} \quad and \quad \frac{1}{q}+\frac{1}{p}=\frac{1}{4}.
\end{equation}
\end{corollary}
\noindent The following estimate is more precise and in fact it implies Theorem 3.6 by duality.
\begin{lemma}
For any $h\in C_c^{\infty}(\mathbb{R}_x^2\times \mathbb{T}_y\times \mathbb{R}_t)$, there holds that 
\begin{equation}\label{equation}
\aligned
&\|\int_{s\in \mathbb{R}} e^{-is\Delta_{\mathbb{R}^2\times \mathbb{T}}} P_{\leq N} h(x,y,s) ds\|_{L_{x,y}^{2}(\mathbb{R}_x^2\times \mathbb{T}_y\times \mathbb{R}_t)} \\ &\lesssim N^{\frac{3}{2}-\frac{5}{p}} \|h\|_{l_{\gamma}^2 L_{x,y,t}^{p^{'}}(\mathbb{R}^2 \times \mathbb{T} \times [2\pi\gamma,2\pi(\gamma+1)])}+N^{1-\frac{16}{5p}} \|h\|_{l_{\gamma}^{q^{'}} L_{x,y,t}^{p^{'}}(\mathbb{R}^2\times \mathbb{T} \times [2\pi\gamma,2\pi(\gamma+1)])}  
\endaligned
\end{equation}
for any (p,q) satisfies (3.2).
\end{lemma}
\noindent The proof of Lemma 3.8 is similar to [14, Lemma 3.2] and Lemma 3.2 (see also [38, Lemma 3.2]). Since the spirits of the proofs are similar, we omit it.\vspace{3mm}
\section{Local theory and small-data result}
\noindent \textbf{For cubic $\mathbb{R}^3\times \mathbb{T}$ problem:} \vspace{3mm}

\noindent We define ``$Z$-norm" (scattering norm) as follows
\[ 
\|u\|_{Z(I)}=(\sum_{N\geq 1}N^{6-p_0} \|1_I(t) P_N u\|_{l^{\frac{2p_0}{p_0-3}}L_{x,y,t}^{p_0}(\mathbb{R}^3\times \mathbb{T} \times I_{\gamma})}^{p_0})^\frac{1}{p_0}.
\]
\noindent Here $p_0$ is a constant satisfying $5<p_0<\frac{11}{2}$. It is not hard to verify that the scattering norm of the linear solution is control by the $H^1$-norm of the initial data by using Stricharz estimate (Corollary 3.3). For convenience, we define ``$Z^{'}$-norm" which is a mixture of $Z$-norm and $X^{1}$-norm as follows
\begin{equation}\label{equation}
\|u\|_{Z^{'}(I)}=\|u\|_{Z(I)}^{\frac{3}{4}} \|u\|_{X^1(I)}^{\frac{1}{4}}.
\end{equation}
\noindent Now we are ready to prove the local well-posedness and small-data scattering of (1.2). 
\begin{lemma}[Bilinear Estimate] Suppose that $ u_i=P_{N_i}u$, for $i=1,2$ satisfying $N_1\geq N_2$. There exists $\delta>0$ such that the following estimate holds for any interval $I\in \mathbb{R}$:   
\begin{equation}\label{equation}
\|u_1 u_2\|_{L_{x,t}^2(\mathbb{R}^3\times \mathbb{T}\times I)} \lesssim (\frac{N_2}{N_1}+\frac{1}{N_2})^\delta \|u_1\|_{Y^0(I)} \|u_2\|_{Z^{'}(I)} .
\end{equation}
\end{lemma}
\noindent \emph{Proof:} Without loss of generality, we can assume that $I=\mathbb{R}$. On one hand, we need the following estimate which follows as in [17, Proposition 2.8],
\begin{equation}\label{equation}
\|u_1 u_2\|_{L^2(\mathbb{R}^3\times \mathbb{T}\times \mathbb{R})} \lesssim N_2(\frac{N_2}{N_1}+\frac{1}{N_2})^\delta \|u_1\|_{Y^{0}(\mathbb{R})} \|u_2\|_{Y^{0}(\mathbb{R})} .
\end{equation}  
\noindent And it suffices to prove the following estimate, if it is hold then we can just combine the two inequalities (noticing the definition of $Z^{'}$-norm) and we will get the lemma completed.
\begin{equation}\label{equation}
\|u_1 u_2\|_{L^2(\mathbb{R}^3\times \mathbb{T} \times \mathbb{R})} \lesssim  \|u_1\|_{Y^{0}(\mathbb{R})} \|u_2\|_{Z(\mathbb{R})} .
\end{equation}  
\noindent We first notice that, by orthogonality considerations, we may replace $u_1$ by $P_C u_1$ where $C$ is a cube of dimension $N_2$. By using H$\ddot{o}$lder's inequality, we have,
\begin{align*}
\|(P_{C}u_{1})u_2\|_{L^2_{x,y,t}} &\lesssim \|P_{C}u_{1}\|_{l_{\gamma}^{\frac{2p_0}{3}}L_{x,y,t}^{\frac{2p_0}{p_0-2}}(\mathbb{R}^3\times \mathbb{T} \times I_{\gamma})} \|u_2\|_{l_{\gamma}^{\frac{2p_0}{p_0-3}}L_{x,y,t}^{p_0}(\mathbb{R}^3\times \mathbb{T}\times I_{\gamma})}      \\
&\lesssim N_2^{\frac{6-p_0}{p_0}}\|P_{C}u_{1}\|_{U_{\Delta}^{\frac{2p_0}{p_0-2}}} \|u_2\|_{l_{\gamma}^{\frac{2p_0}{p_0-3}}L_{x,y,t}^{p_0}(\mathbb{R}^3\times \mathbb{T}\times I_{\gamma})}  \\
&\lesssim \|P_{C}u_{1}\|_{U_{\Delta}^{\frac{2p_0}{p_0-2}}} N_2^{\frac{6-p_0}{p_0}} \|u_2\|_{l_{\gamma}^{\frac{2p_0}{p_0-3}}L_{x,y,t}^{p_0}(\mathbb{R}^3\times \mathbb{T}\times I_{\gamma})}  \\
&\lesssim \|P_{C}u_1\|_{Y^0} \|u_2\|_{Z(\mathbb{R})}.
\end{align*}

\noindent \emph{Remark 1.} The index threshold condition in Stricharz estimate $\frac{22}{7}<\frac{2p_0}{p_0-2}<6$ requires $3<p_0<\frac{11}{2}$. That is a reason why we have the restriction for $p_0$.\vspace{3mm}

\noindent \emph{Remark 2.} We have used another version of Stricharz estimate as follows: for $p>\frac{22}{7}$ and $q$ as in Theorem 3.1, the following estimate holds for any time interval $I\subset \mathbb{R}$ and every cube $Q\subset \mathbb{R}^4$ of size $N$:
\begin{equation}\label{equation}
\|1_{I}(t) P_{Q}u\|_{l_{\gamma}^q L_{x,y,t}^p} \lesssim N^{2-\frac{6}{p}} \|u\|_{U_{\Delta}^{min(p,q)}}(I;L^2(\mathbb{R}^3\times \mathbb{T})). 
\end{equation}

\noindent \emph{Remark 3.} By using Stricharz estimate and the embedding properties of the function spaces, it is not hard to verify $Z$-norm is weaker than $X^1$-norm.\vspace{3mm}

\noindent Based on bilinear estimate, we can prove the following nonlinear estimate which is a crucial step of the local theory. The proof of the following lemmas and theorems follows from standard arguments as in [14, 20, 38].
\begin{lemma} [Nonlinear Estimate] For $u_i \in X^1(I)$, $i=1,2,3$. There holds that
\begin{equation}\label{equation}
\|\widetilde{u}_1 \widetilde{u}_2 \widetilde{u}_3 \|_{N(I)} \leq \sum_{(i,j,k)=(1,2,3) } \|u_i\|_{X^{1}(I)} \|u_j\|_{Z^{'}(I)} \|u_k\|_{Z^{'}(I)} 
\end{equation} where $\tilde{u}_i$ is either $u_i$ or $\bar{u}_i$.
\end{lemma}
\noindent \emph{Proof:} It suffices to prove the following estimate: (Without loss of generality, let $I = \mathbb{R}$)
\begin{equation}\label{equation}
\|\sum_{K\geq 1}P_{K} u_1 \prod_{i=2,3}P_{\leq CK} \tilde{u}_{i}\|_{N(\mathbb{R})}\lesssim_C \|u_1\|_{X^{1}(\mathbb{R})} \|u_2\|_{Z^{'}(\mathbb{R})} \|u_3\|_{Z^{'}(\mathbb{R})} .
\end{equation}
\noindent Using Theorem 2.1, it suffices to prove for any $u_0\in Y^{-1}$ and $||u_0||_{Y^{-1}}\leq 1$ 
\begin{equation}\label{equation}
\sum_{N_{1}} |\int_{\mathbb{R}^3\times \mathbb{T}\times \mathbb{R}} \bar{u}_0 P_{N_1}u_1 \prod_{i=2,3}(P_{\leq CN_{1}} \tilde{u}_{i}) dx dy dt| \leq \|u_0\|_{Y^{-1}} \|u_1\|_{X^1(\mathbb{R})} \|u_2\|_{Z^{'}(\mathbb{R})} \|u_3\|_{Z^{'}(\mathbb{R})}  .
\end{equation}
\noindent Now we split them as follows, let $u_i=\sum\limits_{N_i\geq 1}P_{N_i} u_i$, $i=0,1,2,3$, denoting $u_j^{N_j}= P_{N_j} u_j$ and then the estimate would follow from the following bound:
\begin{equation}\label{equation}
\sum_{S(N_0,N_1,N_2,N_3)} |\int u_0^{N_0} u_1^{N_1} u_2^{N_2} u_3^{N_3}dxdydt| \lesssim \|u_0\|_{Y^{-1}}\|u_1\|_{X^{1}}\|u_2\|_{Z^{'}}\|u_3\|_{Z^{'}} .
\end{equation}
\noindent Here we have set of index $S$ to be $\{(N_0,N_1,N_2,N_3) :  N_1 \sim max(N_2,N_0)\geq N_2\geq N_3 \}$ and we split $S$ into the disjoint union of $S_1$ and $S_2$ and $S_1$ is for the elements in $S$ that satisfy $N_1 \sim N_0$ and $S_2$ if for the elements in $S$ that satisfy $N_1 \sim N_2$. And we will estimate $S_1$ and $S_2$ separately. We omit the proof for $S_2$ part since the estimate is similar.\vspace{3mm}

\noindent By using bilinear estimate (4.2) and some basic inequalities and the properties of function spaces, we have, for a term in $S_1$:
\begin{align*}
|\int u_0^{N_0} u_1^{N_1} u_2^{N_2} u_3^{N_3}dxdydt|&\leq \|u_0^{N_0} u_2^{N_2}\|_{L^2} \|u_1^{N_1} u_3^{N_3}\|_{L^2} \\
&\leq (\frac{N_2}{N_0}+\frac{1}{N_2})^{\delta} (\frac{N_3}{N_1}+\frac{1}{N_3})^{\delta} \|u_0^{N_0}\|_{Y^{0}(\mathbb{R})} \|u_1^{N_1}\|_{Y^{0}(\mathbb{R})} \|u_2^{N_2}\|_{Z^{'}(\mathbb{R})} \|u_3^{N_3}\|_{Z^{'}(\mathbb{R})}.
\end{align*}
\noindent By using Cauchy-Schwarz inequality, the sum of the terms in $S_1$:
\begin{align*}
S_1&\lesssim \sum_{N_1\sim N_0} (\frac{N_2}{N_0}+\frac{1}{N_2})^{\delta} (\frac{N_3}{N_1}+\frac{1}{N_3})^{\delta} \|u_0^{N_0}\|_{Y^{0}(\mathbb{R})} \|u_1^{N_1}\|_{Y^{0}(\mathbb{R})} \|u_2^{N_2}\|_{Z^{'}(\mathbb{R})} \|u_3^{N_3}\|_{Z^{'}(\mathbb{R})} \\
&\lesssim (\sum_{N_1\sim N_0} \frac{N_0}{N_1} \|u_0^{N_0}\|_{Y^{-1}(\mathbb{R})} \|u_1^{N_1}\|_{Y^{1}(\mathbb{R})}) \|u_2\|_{Z^{'}(\mathbb{R})} \|u_3\|_{Z^{'}(\mathbb{R})} \\
&\lesssim \|u_0\|_{Y^{-1}(\mathbb{R})} \|u_1\|_{X^{1}(\mathbb{R})} \|u_2\|_{Z^{'}(\mathbb{R})} \|u_3\|_{Z^{'}(\mathbb{R})}.
\end{align*}
\noindent This finishes Lemma 4.2.
\begin{theorem}\label{}[Local Well-posedness] Let $E > 0$ and $\|u_0\|_{ H^1(\mathbb{R}^3 \times \mathbb{T}) }<E$,
then there exists $\delta_0=\delta_0(E)>0$ such that if 
\[\|e^{it\Delta}u_0\|_{Z(I)}< \delta \] 
\noindent for some $ \delta \leq \delta_0 $, $0\in I$. Then there exists a unique strong solution $u\in X_{c}^{1}(I)$ satisfying $u(0)=u_0$ and we can get an estimate,
\begin{equation}\label{equation}
\|u(t)-e^{it\Delta_{\mathbb{R}^3\times \mathbb{T}}}u_0\|_{X^1(I)}\leq (E\delta)^{\frac{3}{2}} .
\end{equation}  
\end{theorem}
\noindent \emph{Remark 1.} Observe that if $u \in X^1_c(\mathbb{R})$, then $u$ scatters as $t \rightarrow \pm\infty $ as in (1.3). Also, if $E$ is small enough, $I$ can be taken to $\mathbb{R}$ which proves the small data scattering of (1.1).\vspace{3mm}

\noindent \emph{Proof:} First, we consider a mapping defined as follows,
\[
\Phi(u)=e^{it\Delta} u_0-\int_{0}^{t}{e^{i(t-s)\Delta} |u(s)|^2 u(s)ds}.
\]
\noindent And we define a set $ B=\{u\in X_{c}^{1}(I) :   \|u\|_{X^1(I)}\leq 2E  \textmd{ and }  \|u\|_{Z(I)}\leq 2\delta \}$. Now we will verify two properties of $\Phi$: 1. $\Phi$ maps $B$ to $B$. 2. $\Phi$ is a contraction mapping. \vspace{3mm}

\noindent 1. For $ u \in B$, we can use the nonlinear estimate in Lemma 4.2 and let $\delta \leq 1$ and small enough to make $E^3 \delta$ small enough, we have:                
\[
\|\Phi(u)\|_{X^1(I)} \leq \|e^{it\Delta} u_0\|_{X^1(I)}+\||u|^2 u\|_{N(I)} \leq E+CE^{\frac{3}{2}} \delta^{\frac{3}{2}} \leq 2E,
\]
\[
\|\Phi(u)\|_{Z(I)} \leq \|e^{it\Delta} u_0\|_{Z(I)}+\||u|^2 u\|_{N(I)} \leq \delta+CE^{\frac{3}{2}}\delta^{\frac{3}{2}} \leq 2\delta.
\]
\noindent 2.\begin{align*}
\|\Phi(u)-\Phi(v)\|_{X^{1}(I)} &\lesssim \|u-v\|_{X^1(I)} (\|u\|_{X^1(I)}+\|v\|_{X^1(I)})
(\|u\|_{Z^{'}(I)}+\|v\|_{Z^{'}(I)})    \\
&\leq  C\|u-v\|_{X^1(I)} E^{\frac{5}{4}} \delta^{\frac{3}{4}}  \\
&\leq  C\|u-v\|_{X^1(I)} [(E^3\delta)^\frac{5}{12} \delta^\frac{1}{3}] \\
&\leq C\frac{1}{2}\|u-v\|_{X^1(I)}.
\end{align*}

\noindent Thus the result now follows from the Picard's fixed point argument.
\begin{theorem}\label{theorem}
[Controlling Norm] Let $u\in X_{c,loc}^{1}(I)$ be a strong solution on $I\in \mathbb{R} $ satisfying 
\begin{equation}\label{equation}
\|u\|_{Z(I)}< \infty .
\end{equation}
\noindent Then we have two conclusions, \vspace{3mm}

\noindent (1) If $I$ is finite, then u can be extended as a strong solution in $X_{c,loc}^{1}(I^{'})$ on a strictly larger interval $I^{'}$ ,$I\subsetneq I^{'}\subset \mathbb{R}$. In particular, if u blows up in finite time, then the $Z$ norm of $u$ has to blow up.  \vspace{3mm}

\noindent (2) If $I$ is infinite, then $ u\in X_{c}^{1}(I)$.     
\end{theorem}
\noindent \emph{Proof:} Without loss of generality, for the finite case we can assume $I=[0,T)$ and we want to extend it to $[0,T+v)$ for some $v>0$. Denoting $E=\sup\limits_{I}\|u(t)\|_{H^{1}(\mathbb{R}^3 \times \mathbb{T})}$ and using the time-divisibility of `$Z$-norm', there exists $T_1$ such that $T-1<T_1<T$ such that
\[\|u\|_{Z([T_{1},T])} \leq \epsilon,\]
\noindent where $\epsilon$ is to be decided. This allows to conclude:
\[\|u(t)-e^{i(t-T_{1})\Delta } u(T_1) \|_{X^1([T_{1},T])} \lesssim \|u\|_{X^1([T_1,T))}^{\frac{3}{2}} \|u\|_{Z([T_1,T))}^{\frac{3}{2}}\leq C\epsilon^{\frac{3}{2}}  \|u\|_{X^1([T_1,T))}^{\frac{3}{2}}.\]
\noindent By bootstrap argument, we get,

\[ \|u\|_{X^1([T_1,T))}^{\frac{3}{2}} \lesssim E.  \]

\noindent If $\epsilon$ is small enough and, making $\epsilon$ possibly smaller, we have,
\begin{align*}
\|e^{i(t-T_1)\Delta} u(T_1)\|_{Z([T_1,T))} &\leq \|u\|_{Z([T_1,T))} +\|e^{i(t-T_1)\Delta} u(T_1)-u(t)\|_{Z([T_1,T))}   \\
&\leq \epsilon+\|e^{i(t-T_1)\Delta} u(T_1)-u(t)\|_{X^{1}([T_1,T))}   \\
&\leq \epsilon+C^{'}\epsilon^{\frac{3}{2}} E^{\frac{3}{2}}    \\
&\leq \frac{3}{4} \delta_{0}(E).
\end{align*}

\noindent Notice that we can let $\epsilon$ small enough s.t. $\epsilon<\frac{1}{4}$ and $\epsilon E<(\frac{1}{2} \delta_{0}(E))^{\frac{2}{3}}$. This allows to find an interval $[T_1,T+v]$ for which :
\[ \|e^{i(t-T_1)\Delta}u(T_1)\|_{Z([T_1,T+v])}<\delta_0 .  \]
\noindent That finishes the proof by using the Theorem 4.3. Moreover, by using the symmetries of the equation, the above argument also covers the case when $I$ is an arbitrary bounded interval.\vspace{3mm}

\noindent Now we turn to the infinite case. Without loss of generality, it is enough to consider the case $I=(a,\infty)$. Choosing $T$ to be large enough so that
\[||u||_{Z([T,\infty))} \leq \epsilon,
\]
\noindent we get that for any $T^{'}>T$: 
\[\|u(t)-e^{i(t-T)\Delta } u(T) \|_{X^1([T,T^{'}))} \lesssim \|u\|_{X^1([T,T^{'}))}^{\frac{3}{2}} \|u\|_{Z([T_1,T))}^{\frac{3}{2}}\leq C\epsilon^{\frac{3}{2}}  \|u\|_{X^1([T,T^{'}))}^{\frac{3}{2}},\]

\noindent which gives that $||u||_{X^1([T,T^{'}))} \lesssim E$ for any $T^{'}>T$ and we have
\[||e^{i(t-T)\Delta}u(T)||_{Z([T,\infty))} \leq 2\epsilon \leq \delta_0(E)
\]
\noindent if $\epsilon$ small enough. The result now follows from by using Theorem 4.3. 
\begin{theorem}\label{theorem}[Stability Theory] Let $ I\in \mathbb{R}$ be an interval, and let $\widetilde{u} \in X^1(I) $ solve the approximate solution,
\begin{equation}\label{equation}
(i\partial_{t}+\Delta_{\mathbb{R}^3 \times \mathbb{T}} )\widetilde{u}=\rho |\widetilde{u}|^2\widetilde{u}+e \textmd{ and }  \rho\in [0,1]. 
\end{equation}
\noindent Assume that:
\begin{equation}\label{equation}
\|\widetilde{u}\|_{Z(I)}+\|\widetilde{u}\|_{L_{t}^{\infty}(I,H^{1}(\mathbb{R}^3\times \mathbb{T}))} \leq M.
\end{equation}
\noindent There exists $\epsilon_{0}=\epsilon_{0}(M) \in (0,1] $ such that if for some $t_{0} \in I$:
\begin{equation}\label{equation}
\|\widetilde{u}(t_{0})-u_{0}\|_{H^{1}(\mathbb{R}^3\times \mathbb{T})}+\|e\|_{N(I)} \leq \epsilon < \epsilon_{0},
\end{equation}
\noindent then there exists a solution $u(t)$ to the exact equation:
\begin{equation}\label{equation}
(i\partial_{t}+\Delta_{\mathbb{R}^3 \times \mathbb{T}} )u=|u|^2u \end{equation}
\noindent with initial data $u_0$ satisfies 
\begin{equation}\label{equation} \|u\|_{X^1(I)}+\|\widetilde{u}\|_{X^1(I)} \leq C(M),   \quad  \quad  \|u-\widetilde{u}\|_{X^1(I)} \leq C(M)\epsilon. \end{equation}
\end{theorem}
\noindent \emph{Proof:} The proof is very similar to the proof of [14, Proposition 4.7]. The proof relies tightly on the estimate of Lemma 4.2 (nonlinear estimate) and the trick of division of the intervals.\vspace{3mm} 

\noindent First, we consider for an interval $J\in I$ s.t. $\|\widetilde{u}\|_{Z(J)} \leq \epsilon $ (That is the additional smallness assumption, $\epsilon$ is to be decided). We will prove the theorem under this assumption. By local existing argument for the approximate equation, there exists $ \delta_{1}(M)$ that if 
\[ \|e^{i(t-t_{*})\Delta} \widetilde{u}(t_{*})\|_{Z(J)}+\|e\|_{N(J)} \leq \delta_{1} \]
for some $t_{*}\in J$, then $\widetilde{u}\in X^1(J)$ is unique and satisfies:
\[ \|\widetilde{u}- e^{i(t-t_{*})\Delta} \widetilde{u}(t_{*})\|_{X^1(J)} \leq C\|\widetilde{u}\|_{X^1(J)}^{\frac{3}{2}} \|\widetilde{u}\|_{Z(J)}^{\frac{3}{2}}+\|e\|_{N(J)}. \]

\noindent We can conclude
\[ \|\widetilde{u}\|_{X^1(J)} \lesssim M+1\quad \textmd{and} \quad \|e^{i(t-t_{*})\Delta} \widetilde{u}(t_{*})\|_{Z(J)} \lesssim \epsilon.\]\noindent if $\epsilon<\epsilon_1(M)$ is small enough.\vspace{3mm}

\noindent Second, let us estimate the difference of the solutions. Consider solution $u$ with initial data $u_{*}$ satisfying $\|u_{*}-\widetilde{u}(t_{*})\|_{H^{1}} \leq \epsilon $ and living on an interval $J_{u} \in J$ containing $t_{*}$. And, we want to prove the following estimate for some constant $C$ independent of $J_u$ to be specified later:
\begin{equation}\label{equation}
\|u-\widetilde{u} \|_{X^1(J_u)}\leq C\epsilon. 
\end{equation}

\noindent Let $w=u-\widetilde{u} $, then we know that $w$ satisfies: 
\[
(i\partial_t +\Delta)w=\rho(|\widetilde{u}+w|^2 (\widetilde{u}+w)-|\widetilde{u}|^2 \widetilde{u} )-e.
\]
\noindent Adopting the bootstrap hypothesis: 
\[ \|w\|_{X^1(J_u\cap[t_{*}-t,t_{*}+t]) }\leq2C\epsilon. \]
\noindent For convenience, we denote $J_u\cap[t_{*}-t,t_{*}+t) $ by $J_t$, by using nonlinear estimate, we compute:
\begin{align*}
\|w\|_{X^1(J_t)}&\lesssim \|u(t_{*})-\widetilde{u}(t_{*})\|_{H^1(\mathbb{R}^3\times \mathbb{T})}+\|w\|_{X^1(J_t)} \| \widetilde{u}\|_{X^1(J_t)}  \| \widetilde{u}\|_{Z^{'}(J_t)}+\|e\|_{N(J_t)}     \\
&\lesssim \epsilon+\|w\|_{X^1(J_t)} \|\widetilde{u}\|_{X^1(J_t)}^{\frac{5}{4}}   \|\widetilde{u}\|_{Z(J_t)}^{\frac{3}{4}}   \\
&\leq C_1 \epsilon+C_1 M^{\frac{5}{4}} \epsilon^{\frac{3}{4}} \|w\|_{X^1(J_t)}.     
\end{align*}
\noindent As a result, if $\epsilon<\epsilon_1(M)$ with $\epsilon_1(M)$ small enough in terms of $M$, we conclude that $||u-\tilde{u}||_{X^1(J_t)}\leq 2C_1\epsilon$, which close the the bootstrap argument with $C=2C_1$. This finishes the proof under the smallness assumption. \vspace{3mm}

\noindent Now, to generalize the argument to the whole interval $I$, we split $I$ into $N=C(M,\epsilon_1(M))$ intervals $I_k=[T_k,T_{k+1})$ such that:
\[||u||_{Z(I_k)}\leq \frac{\epsilon_1(M)}{100} \quad \textmd{and} \quad ||e||_{N(I_k)}\leq \frac{\epsilon_1(M)}{100}.
\]
\noindent If $\epsilon_0(M)$ is chosen sufficiently small in terms of $N$, $M$ and $\epsilon_1(M)$, we can iterate the first part of the proof on each interval $I_k$ while keeping the condition
\[||u(T_k)-\tilde{u}(T_k)||_{H^1(\mathbb{R}^3\times \mathbb{T})}+||e||_{N(I_k)}+||u||_{Z(I_k)}<\epsilon_1(M)
\]
\noindent always satisfied for each $k$. This finishes the proof of Theorem 4.5.\vspace{3mm}

\noindent \textbf{For quintic $\mathbb{R}^2\times \mathbb{T}$ problem:} we have the following statements for the quintic analogues:\vspace{3mm}

\noindent We define the scattering norm for this problem as follows:
\begin{equation}
\|u\|_{Z(I)}=\sum_{p_0=8,10}(\sum_{N\geq 1}N^{5-\frac{p_0}{2} } \|1_I(t) P_N u\|_{l^{\frac{4p_0}{p_0-4}}L_{x,y,t}^{p_0}(\mathbb{R}^2\times \mathbb{T} \times I_{\gamma})}^{p_0})^\frac{1}{p_0}.
\end{equation}

\noindent For convenience, we define $Z^{'}$-norm as a mixture of $Z$-norm and $X^1$-norm as follows:
\begin{equation}\label{equation}
\|u\|_{Z^{'}(I)}=\|u\|_{Z(I)}^{\frac{3}{4}} \|u\|_{X^1(I)}^{\frac{1}{4}}.
\end{equation}
\noindent Now we are ready to prove the local well-posedness and small-data scattering of (1.3). We will prove the following trilinear estimate first:
\begin{lemma}[Trilinear estimate]Suppose that $ u_i=P_{N_i}u$, for $i=1,2,3$ satisfying $N_1\geq N_2 \geq N_3$. There exists $\delta>0$ such that the following estimate holds for any interval $I\in \mathbb{R}$:   
\begin{equation}\label{equation}
\|u_1 u_2 u_3\|_{L_{x,t}^2(\mathbb{R}^2\times \mathbb{T}\times I)} \lesssim (\frac{N_3}{N_1}+\frac{1}{N_2})^\delta \|u_1\|_{Y^0(I)} \|u_2\|_{Z^{'}(I)} \|u_3\|_{Z^{'}(I)} .
\end{equation}
\end{lemma}
\noindent \emph{Proof:} Without loss of generality, we can assume that $I=\mathbb{R}$. On one hand, we need the following estimate which follows as in [17, Proposition 2.8],
\begin{equation}\label{equation}
\|u_1 u_2 u_3\|_{L^2(\mathbb{R}^2\times \mathbb{T}\times \mathbb{R})} \lesssim N_2 N_3(\frac{N_3}{N_1}+\frac{1}{N_2})^\delta \|u_1\|_{Y^{0}(\mathbb{R})} \|u_2\|_{Y^{0}(\mathbb{R})} \|u_3\|_{Y^{0}(\mathbb{R})} .
\end{equation}  
\noindent And it suffices to prove the following estimate, if it is hold then we can just combine those two inequalities to get the lemma proved.
\begin{equation}\label{equation}
\|u_1 u_2 u_3\|_{L^2(\mathbb{R}^2\times \mathbb{T} \times \mathbb{R})} \lesssim  \|u_1\|_{Y^{0}(\mathbb{R})} \|u_2\|_{Z(\mathbb{R})} \|u_3\|_{Z(\mathbb{R})} .
\end{equation}  
\noindent By orthogonality considerations, we may replace $u_1$ by $P_C u_1$ where $C$ is a cube of dimension $N_2$. By using H$\ddot{o}$lder's inequality, we have,
\begin{align*}
\|(P_{C}u_{1})u_2 u_3\|_{L^2_{x,y,t}} &\lesssim \|P_{C}u_{1}\|_{l_{\gamma}^{\frac{40}{9}}L_{x,y,t}^{\frac{40}{11}}(\mathbb{R}^2\times \mathbb{T} \times I_{\gamma})} \|u_2\|_{l_{\gamma}^{8}L_{x,y,t}^{8}(\mathbb{R}^2\times \mathbb{T}\times I_{\gamma})} \|u_3\|_{l_{\gamma}^{\frac{20}{3}}L_{x,y,t}^{10}(\mathbb{R}^2\times \mathbb{T}\times I_{\gamma})}     \\
&\lesssim N_2^{\frac{1}{8}}\|P_{C}u_{1}\|_{U_{\Delta}^{\frac{40}{11}}} \|u_2\|_{l_{\gamma}^{8}L_{x,y,t}^{8}(\mathbb{R}^2\times \mathbb{T}\times I_{\gamma})} \|u_3\|_{Z(\mathbb{R})} \\
&\lesssim \|P_{C}u_{1}\|_{U_{\Delta}^{\frac{40}{11}}} N_2^{\frac{1}{8}} \|u_2\|_{l_{\gamma}^{8}L_{x,y,t}^{8}(\mathbb{R}^2\times \mathbb{T}\times I_{\gamma})}\|u_3\|_{Z(\mathbb{R})}  \\
&\lesssim \|P_{C}u_1\|_{Y^0} \|u_2\|_{Z(\mathbb{R})}\|u_3\|_{Z(\mathbb{R})}.
\end{align*}
\noindent This finishes the proof of (4.20).\vspace{3mm}

\noindent Based on the trilinear estimate, we are ready to prove Nonlinear Estimate Lemma, Controlling Norm Theorem and Stability Theorem as follows. We only state those propositions and omit the proofs since the proofs of those propositions are analogues of [14, Lemma 4.3, Proposition 4.5, Proposition 4.6 and Proposition 4.7] once we have the trilinear estimate. (see also [20, 36] and the cubic $\mathbb{R}^3\times \mathbb{T}$ case in this paper)

\begin{lemma} [Nonlinear Estimate] For $u_i \in X^1(I)$, $i=1,2,3,4,5$. There holds that
\begin{equation}\label{equation}
\|\widetilde{u}_1 \widetilde{u}_2 \widetilde{u}_3 \widetilde{u}_4 \widetilde{u}_5\|_{N(I)} \leq \sum_{\{i,j,p,m,n\}=\{1,2,3,4,5\} } \|u_i\|_{X^{1}(I)} \|u_j\|_{Z^{'}(I)} \|u_p\|_{Z^{'}(I)} \|u_m\|_{Z^{'}(I)}\|u_n\|_{Z^{'}(I)}
\end{equation} where $\tilde{u}_i$ is either $u_i$ or $\bar{u}_i$.
\end{lemma}

\begin{theorem}\label{}[Local Well-posedness] Let $E > 0$ and $\|u_0\|_{ H^1(\mathbb{R}^2 \times \mathbb{T}) }<E$,
then there exists $\delta_0=\delta_0(E)>0$ such that if 
\[\|e^{it\Delta}u_0\|_{Z(I)}< \delta \] 
\noindent for some $ \delta \leq \delta_0 $, $0\in I$. Then there exists a unique strong solution $u\in X_{c}^{1}(I)$ satisfying $u(0)=u_0$ and we can get an estimate,
\begin{equation}\label{equation}
\|u(t)-e^{it\Delta_{\mathbb{R}^2\times \mathbb{T}}}u_0\|_{X^1(I)}\leq E^2 \delta^{3} .
\end{equation}  
\end{theorem}

\begin{theorem}\label{theorem}
[Controlling Norm] Let $u\in X_{c,loc}^{1}(I)$ be a strong solution on $I\in \mathbb{R} $ satisfying 
\begin{equation}\label{equation}
\|u\|_{Z(I)}< \infty .
\end{equation}
\noindent Then we have two conclusions, \vspace{3mm}

\noindent (1) If $I$ is finite, then u can be extended as a strong solution in $X_{c,loc}^{1}(I^{'})$ on a strictly larger interval $I^{'}$ ,$I\subsetneq I^{'}\subset \mathbb{R}$. In particular, if u blows up in finite time, then the $Z$-norm of $u$ has to blow up.  \vspace{3mm}

\noindent (2) If $I$ is infinite, then $ u\in X_{c}^{1}(I)$.     
\end{theorem}

\begin{theorem}\label{theorem}[Stability Theory] Let $ I\in \mathbb{R}$ be an interval, and let $\widetilde{u} \in X^1(I) $ solve the approximate solution,
\begin{equation}\label{equation}
(i\partial_{t}+\Delta_{\mathbb{R}^2 \times \mathbb{T}} )\widetilde{u}=\rho |\widetilde{u}|^4\widetilde{u}+e \textmd{ and }  \rho\in [0,1]. 
\end{equation}
\noindent Assume that:
\begin{equation}\label{equation}
\|\widetilde{u}\|_{Z(I)}+\|\widetilde{u}\|_{L_{t}^{\infty}(I,H^{1}(\mathbb{R}^2\times \mathbb{T}))} \leq M.
\end{equation}
\noindent There exists $\epsilon_{0}=\epsilon_{0}(M) \in (0,1] $ such that if for some $t_{0} \in I$:
\begin{equation}\label{equation}
\|\widetilde{u}(t_{0})-u_{0}\|_{H^{1}(\mathbb{R}^2\times \mathbb{T})}+\|e\|_{N(I)} \leq \epsilon < \epsilon_{0},
\end{equation}
\noindent then there exists a solution $u(t)$ to the exact equation:
\begin{equation}\label{equation}
(i\partial_{t}+\Delta_{\mathbb{R}^2 \times \mathbb{T}} )u=|u|^4u \end{equation}
\noindent with initial data $u_0$ satisfies 
\begin{equation}\label{equation} \|u\|_{X^1(I)}+\|\widetilde{u}\|_{X^1(I)} \leq C(M),   \quad  \quad  \|u-\widetilde{u}\|_{X^1(I)} \leq C(M)\epsilon. \end{equation}
\end{theorem}

\section{Nonlinear Analysis of the Profiles}
\noindent In this section, we describe and analyze the profiles that appear in the linear and nonlinear profile decomposition.\vspace{3mm} 

\noindent \textbf{For cubic $\mathbb{R}^3\times \mathbb{T}$ problem:} \vspace{3mm}

\noindent We recall the motivation discussed in Section 1: In view of the scaling-invariant of the IVP (1.1) under 
\[\mathbb{R}_x^3 \times \mathbb{T}_y \rightarrow M_{\lambda}:=\mathbb{R}_x^3 \times (\lambda^{-1}\mathbb{T})_{y},\quad u \rightarrow \tilde{u}(x,y,t)=\lambda u(\lambda x,\lambda y,\lambda^2 t).\]

\noindent When $\lambda \rightarrow 0$, the manifolds $M_\lambda$ will be similar to $\mathbb{R}^4$. The appearance is a manifestation of the energy-critical nature of the nonlinearity. This extreme behavior corresponds to Euclidean profile. Precise description is as follows.\vspace{3mm}

\noindent \emph{Remark.} We also refer to [14, Section 5], [20, Section 4] and [36, Section 5] for more information. For those problems, Euclidean profiles also appear in the analysis of profile decomposition according to the structures of the corresponding equations.\vspace{3mm}

\noindent \textbf{Euclidean Profiles.} The Euclidean profiles define a regime where we can compare solutions of cubic NLS on $\mathbb{R}^4$ with those on $\mathbb{R}^3\times \mathbb{T}$. We fix a spherically symmetric function $\eta \in C_{0}^{\infty}(\mathbb{R}^4)$ supported in the ball of radius 2 and equal to 1 in the ball of radius 1. Given $\phi \in H^1(\mathbb{R}^4)$ and a real number $N \geq 1$, we define:
\[
Q_{N} \phi \in H^1(\mathbb{R}^4)   \quad  \quad  (Q_{N} \phi)  (x)=\eta(\frac{x}{N^{\frac{1}{2}}}) \phi(x),
\]
\begin{equation}\label{euqation}
\phi_{N}\in H^1(\mathbb{R}^4)  \quad  \quad \phi_{N}(x)=N(Q_N \phi)(Nx),
\end{equation}
\[
f_N\in H^1(\mathbb{R}^3 \times \mathbb{T}) \quad \quad f_N(y)=\phi_N(\Psi^{-1}(y)),
\]
\noindent where $\Psi$ is the identity map from the unit ball of $\mathbb{R}^4$ to $\mathbb{R}^3\times \mathbb{T}$. Thus $ Q_N \phi $ is a compactly supported modification of the profile $\phi$, $\phi_{N}$ is an $\dot{H^{1}}$-invariant rescaling of $Q_N \phi $, and $f_{N}$ is the function obtained by transferring $\phi_{N}$ to a neighborhood of $0$ in $\mathbb{R}^3 \times \mathbb{T}$. We notice that 
\[
\|f_{N}\|_{H^1(\mathbb{R}^3 \times \mathbb{T})}  \lesssim \|\phi\|_{\dot{H^{1}}(\mathbb{R}^4)}.
\]
\noindent  And we use scattering result for 4d energy critical NLS by E. Ryckman and M. Visan ([26]) in the following form:
\begin{theorem}\label{theorem}
Assume $\psi \in \dot{H^{1}}(\mathbb{R}^4)$, then there is a unique global solution $v \in C(\mathbb{R}:\dot{H^{1}}(\mathbb{R}^4))$ of the initial-value problem 
\begin{equation}\label{euqation}
(i\partial_{t} +\Delta_{\mathbb{R}^4})v=v|v|^{2}, \quad v(0)=\psi,
\end{equation}
\noindent and
\begin{equation}\label{euqation}
\||\nabla_{\mathbb{R}^4}v|\|_{(L_{t}^{\infty}L_{x}^{2} \cap L_{t}^{2}L_{x}^{4})(\mathbb{R}^4\times \mathbb{R})} \leq \tilde{C}(E_{\mathbb{R}^4}(\psi)).
\end{equation}
\noindent Moreover this solution scatters in the sense that there exists $\psi^{\pm \infty} \in \dot{H^{1}}(\mathbb{R}^4)$ such that 
\begin{equation}\label{euqation}
\|v(t)-e^{it\Delta} \psi^{\pm \infty}\|_{\dot{H^{1}}(\mathbb{R}^4)} \to 0
\end{equation}
as $t \to \pm \infty$. Besides if $\psi \in H^{5}(\mathbb{R}^4)$, then $v \in C(\mathbb{R}:H^5(\mathbb{R}^4))$ and 
\[
\sup\limits_{t\in \mathbb{R}} \|v(t)\|_{H^{5}(\mathbb{R}^4)}  \lesssim_{\|\psi\|_{H^{5}(\mathbb{R}^4)}}   1.
\] \end{theorem}
\noindent Based on the above result, we have:
\begin{theorem}\label{theorem}
\noindent Assume $\phi \in \dot{H^{1}}(\mathbb{R}^4)$, $T_{0} \in (0,\infty)$, and $\rho \in \{0,1\}$ are given, and we define $f_{N}$ as before. Then the following conclusions hold:\vspace{3mm} 

\noindent $(1)$ There is $N_0=N_0 (\phi,T_0)$ sufficiently large such that for any $N\geq N_0$, there is a unique solution $U_{N} \in C((-T_0 N^{-2},T_0 N^{-2}); H^1(\mathbb{R}^3\times \mathbb{T}))$ of the initial-value problem 
\begin{equation}\label{euqation}
(i\partial_t+\Delta)U_{N}=\rho U_{N} |U_{N}|^{2} ,\quad and \quad  U_{N}(0)=f_{N}.
\end{equation}
\noindent Moreover, for any $N \geq N_{0}$, 
\begin{equation}\label{euqation}
\|U_N\|_{X^{1}(-T_0 N^{-2},T_0 N^{-2})} \lesssim_{E_{\mathbb{R}^4}(\phi)} 1.
\end{equation}
$(2)$ Assume $\epsilon_{1} \in (0,1]$ is sufficiently small (depending only on $E_{\mathbb{R}^4}(\phi)$), $\phi ^{'} \in H^{5}(\mathbb{R}^4)$, and $\|\phi-\phi^{'}\|_{\dot{H^{1}}(\mathbb{R}^4)} \leq \epsilon_{1}$. Let $v^{'} \in C(\mathbb{R}:H^5)$ denote the solution of the initial-value problem 
\[
(i\partial_t+\Delta_{\mathbb{R}^4}) v^{'}=\rho v^{'} |v^{'}|^{2}, \quad v^{'} (0)=\phi^{'}.
\] 
\noindent For $R \geq 1$ and $N \geq 10R$, we define 
\[ v_{R}^{'}(x,t)=\eta(\frac{x}{R})v^{'}(x,t)    \quad \quad (x,t) \in \mathbb{R}^{4}\times (-T_0,T_0),
\]
\begin{equation}\label{euqation}
v_{R,N}^{'}(x,t)=N v_{R}^{'}(Nx,N^2 t) \quad \quad (x,t) \in \mathbb{R}^{4}\times (-T_0 N^{-2},T_0 N^{-2}),
\end{equation}
\[
V_{R,N}(y,t)=v_{R,N}^{'}(\Psi^{-1}(y),t) \quad \quad (y,t) \in \mathbb{R}^3\times \mathbb{T} \times (-T_0 N^{-2},T_0 N^{-2}).
\]
\noindent Then there is $R_0 \geq 1$ (depending on $T_0$ and $\phi^{'}$ and $\epsilon_1$) such that, for any $R \geq R_{0}$ and $N \geq 10R$,
\begin{equation}\label{euqation}
\limsup\limits_{N \to \infty} \|U_N-V_{R,N}\|_{X^{1}(-T_0 N^{-2},T_0 N^{-2})}  \lesssim_{E_{\mathbb{R}^4}(\phi)} \epsilon_1.
\end{equation} 
\end{theorem}
\noindent \emph{Proof:} It suffices to prove part (2). All implicit constants are allowed to depend on $\|\phi\|_{\dot{H^{1}}(\mathbb{R}^4)}$. The idea of the proof is to show that with $ R_{0}$ chosen large enough, $V_{R,N}$ is an approximate solution. First, we define: 
\[ e_{R}(x,t) := (i\partial_t +\Delta_{\mathbb{R}^4})v_{R}^{'}-\rho {|v_{R}^{'}|}^{2} v_{R}^{'}.
\]
\noindent Using the fact that $\sup\limits_{t} \|v^{'}(t)\|_{H^5} \lesssim_{\|\phi^{'}\|_{H^5}} 1$, we get that:  
\[
|e_{R}(t,x)|+|\nabla_{\mathbb{R}^4} e_{R}(t,x)| \lesssim 1_{[R/2,4R]} (|v^{'}(t,x)|+|\nabla_{\mathbb{R}^4}v^{'}(t,x)|+|\Delta_{\mathbb{R}^4}v^{'}(t,x)|),
\]
\noindent which directly gives that there exists $R_0 \geq 1$ such that for all $R > R_{0}$   
\[
\lim\limits_{R\to \infty} \| |e_{R}|+|\nabla_{\mathbb{R}^4} e_{R}| \|_{L_{t}^{1} L_{x}^{2}(\mathbb{R}^4 \times (-T,T))} = 0.
\]
\noindent Letting 
\[
e_{R,N}(x,t):=(i\partial_{t}+\Delta_{\mathbb{R}^4})v_{R,N}^{'}-\rho|v_{R,N}^{'}|^{2}v_{R,N}^{'},
\]
\noindent we have that for any $R>R_{0}$ and $N \geq 1$: 
\begin{equation}\label{euqation}
 \| |e_{R,N}|+|\nabla_{\mathbb{R}^4} e_{R,N}| \|_{L_{t}^{1} L_{x}^{2}(\mathbb{R}^4 \times (-TN^{-2},TN^{-2}))} \leq  2\epsilon_{1}
\end{equation}
\noindent with $V_{R,N} $ defined on $\mathbb{R}^3 \times \mathbb{T} \times (-TN^{-2},TN^{-2})$. We let 
\[ E_{R,N}(y,t)=(i\partial_{t}+\Delta_{\mathbb{R}^4})V_{R,N}-\rho|V_{R,N}|^{2}V_{R,N}=e_{R,N}(\Psi^{-1}(y),t). 
\]
\noindent For $R >R_0$ and $N \geq 10R$:
\[  \| |E_{R,N}|+|\nabla_{\mathbb{R}^4} E_{R,N}| \|_{L_{t}^{1} L_{x}^{2}(\mathbb{R}^4 \times (-TN^{-2},TN^{-2}))} \lesssim  \epsilon_{1}
\]
\noindent from which  it follows (using Theorem 2.1) that:
\[ ||E_{R,N}||_{N(-TN^{-2},TN^{-2})} \lesssim \epsilon_{1}.
\]\noindent To verify the requirements of Theorem 4.5, we use (5.3) to conclude that:
\[ ||V_{R,N}||_{L_t^{\infty}H^1(\mathbb{R}^3 \times \mathbb{T} \times(-TN^{-2},TN^{-2}))} \lesssim 1 .
\]
\noindent As for the $Z$-norm control, we choose $N$ to be big enough so that $TN^{-2} \leq \frac{1}{2}$ which makes all summations in the  $Z$-norm consist of at most two terms, after which we estimate the $Z$-norm by using Littlewood-Paley theory and Sobolev embedding theorem as follows: 
\[  ||K^{\frac{6}{p_0}-1}||P_k V_{R,N}||_{L^{p_0}_{x,t}(\mathbb{R}^3\times \mathbb{T} \times (-TN^{-2},TN^{-2}))}||_{l_k^{p_0}} \lesssim ||(1-\Delta)^{\frac{3}{p_0}-\frac{1}{2}}V_{R,N}||_{L^{p_0}_{t,x} } \lesssim  ||(1-\Delta)^{\frac{1}{2}}v^{'}_{R,N}||_{L^{p_0}_{t}L^{\frac{2p_0}{p_0-1}}_{x}} \lesssim_{E(\phi)} 1.
\]

\noindent At last, we know for $R_0$ big enough and $R >R_0$, $N \geq 10R$,
\[ ||f_N-V_{R,N}(0)||_{H^1(\mathbb{R}^3 \times \mathbb{T})} \lesssim ||Q_N\phi-\phi||_{\dot{H}^1(\mathbb{R}^4)}+||\phi^{'}-\phi||_{\dot{H}^1(\mathbb{R}^4)}+||\phi^{'}-V^{'}_{R}(0)||_{\dot{H}^1(\mathbb{R}^4)} \lesssim \epsilon_1.
\]
\noindent This completes the verification of the requirements of Theorem 4.5 which concludes the proof.
\begin{lemma}[Extinction Lemma]Suppose that $\phi \in \dot{H^1}(\mathbb{R}^4), \epsilon > 0$, and  $I \subset \mathbb{R} $ is an interval. Assume that
\begin{equation}\label{euqation}
||\phi||_{\dot{H^1} (\mathbb{R}^4)} \leq 1, \quad ||\nabla_{x} e^{it\Delta}\phi||_{L_t^2 L_x^4(\mathbb{R}^4 \times I)} \leq \epsilon.
\end{equation}
\noindent For $N \geq 1$, we define as before: 
\[Q_N \phi=\eta(N^{-1/2}x)\phi(x),\quad \phi_N=N(Q_N\phi)(Nx),\quad f_N(y)=\phi_N(\Psi^{-1}(y)).
\]
\noindent Then there exists $N_0=N_0(\phi,\epsilon)$ such that for any $N \geq N_0$,
\[ ||e^{it\Delta} f_{N}||_{Z(N^{-2}I)} \lesssim \epsilon.
\]
\end{lemma}
\noindent \emph{Proof:} It suffices to prove that there exists $T_0$ such that for any $N>1$:
\begin{equation}\label{euqation}
||e^{it\Delta} f_N||_{Z(\mathbb{R} \setminus (-N^{-2}T_0,N^{-2}T_0))} \lesssim \epsilon
\end{equation}
\noindent as the rest follows from Lemma 5.2 (with $\rho=0$). Without loss of generality, by limiting arguments, we may assume that $\phi \in C^{\infty}_{0}(\mathbb{R}^4)$. We have, for any $p$, 
\[ f_{N,p}(x)=\frac{1}{2\pi} \int_{\mathbb{T}} \phi_{N}(x,y)e^{-i\langle y,p \rangle} dy= \frac{N}{2\pi} \int_{\mathbb{R}^3} e^{-i\langle y,p \rangle} \phi(Nx,Ny) dy.
\]
\noindent And using dispersive estimate and unitarity, we have 
\begin{equation}\label{euqation} ||e^{it\Delta}P_M f_N(t)||_{L^{\infty}_{x,y}(\mathbb{R}^3\times \mathbb{T})} \lesssim \sup_{x \in \mathbb{R}^3}\sum_{|p| \leq M}|e^{it\Delta_{x}}f_{N,p}(x)| \lesssim \frac{M}{|t|^{\frac{3}{2}}}||f_N||_{L^1_{x,y}} \lesssim \frac{MN^{-3}}{|t|^{\frac{3}{2}}} 
\end{equation}
\noindent and 
\begin{equation}\label{euqation} ||e^{it\Delta}P_M f_N(t)||_{L^{2}_{x,y}(\mathbb{R}^3\times \mathbb{T})}= ||P_M f_N(t)||_{L^{2}_{x,y}(\mathbb{R}^3\times \mathbb{T})} \lesssim M^{-l} ||(1-\Delta)^{\frac{l}{2}} \phi_N||_{L^2(\mathbb{R}^4)} \lesssim M^{-l}N^{l-1}.
\end{equation}
\noindent Then by interpolation we have (choose $l=0, 10000$):
\begin{equation}\label{equation}
\aligned
||e^{it\Delta}P_M f_N(t)||_{L^{p}_{x,y}(\mathbb{R}^3\times \mathbb{T})}  &\lesssim \frac{N^{-2+\frac{4}{p}}}{|t|^{\frac{3}{2}(1-\frac{2}{p})}} [(\frac{M}{N})^{1-\frac{2}{p}-\frac{2l}{p}}] \\
&\lesssim \frac{N^{-2+\frac{4}{p}}}{|t|^{\frac{3}{2}(1-\frac{2}{p})}} min[(\frac{M}{N})^{1-\frac{2}{p}},(\frac{N}{M})^{100}].
\endaligned
\end{equation}
\noindent As a result,
\begin{equation}\label{euqation} 
(\sum_{M} M^{\frac{2}{5}}||e^{it\Delta}P_M f_N||^2_{L_{x,y,t}^{5}(\mathbb{R}^3\times \mathbb{T}\times \{N^2|t|\geq T\})})^{\frac{1}{2}}\lesssim T^{-\frac{4}{5}}. 
\end{equation}
\noindent Also, by using Stricharz estimate, for $p>p_0>5$, we have
\begin{equation}
(\sum_{M} M^{\frac{12}{p}-2}||e^{it\Delta}P_M f_N||^2_{l_{\gamma}^{\frac{2p}{p-3}}L_{x,y,t}^{p}(\mathbb{R}^3\times \mathbb{T} \times \mathbb{R})})^{\frac{1}{2}} \lesssim ||\phi_{N}||_{H^1} \lesssim 1.
\end{equation}
\noindent By interpolation, we can obtain (5.11). That finishes the proof of Lemma 5.3.\vspace{3mm}

\noindent Now we are ready to describe the nonlinear solutions of (1.1) corresponding to data concentrating at a point. Let $\tilde{\mathcal{F}}_e$ denote the set of renormalized Euclidean frames as follows:\vspace{3mm}

\noindent $\tilde{\mathcal{F}}_e := \{(N_k,t_k,x_k)_{k\geq 1}:N_k \in [1,\infty),x_k \in \mathbb{R}^3 \times \mathbb{T} ,N_k \rightarrow \infty$, and either $t_k=0$ for any $k\geq 1$ or $\lim\limits_{k\rightarrow \infty}N_k^2 |t_k|=\infty \}$.\vspace{3mm}

\noindent Given $f \in L^2(\mathbb{R}^3 \times \mathbb{T})$, $t_0 \in \mathbb{R}$, and $x_0 \in \mathbb{R}^3\times \mathbb{T}$, we define:
\begin{equation}\label{euqation} 
\pi_{x_0}f=f(x-x_0), \quad \Pi_{(t_0,x_0)}f=(e^{-it_0\Delta_{\mathbb{R}^3 \times \mathbb{T}}}f)(x-x_0)=\pi_{x_0}e^{it_0\Delta_{\mathbb{R}^3 \times \mathbb{T}}}f. 
\end{equation}
\noindent Also for $\phi \in \dot{H}^1(\mathbb{R}^4)$ and $N \geq 1$, we denote the function obtained in (5.1) by:
\begin{equation}
T^e_{N_k}:=N\tilde{\phi}(N\Psi^{-1}(x)) \quad where \quad \tilde{\phi}(y):=\eta(\frac{y}{N^{\frac{1}{2}}}) \phi(y)
\end{equation}
\noindent and as before observe that $T_N^e:\dot{H}^1(\mathbb{R}^4) \rightarrow H^1(\mathbb{R}^3 \times \mathbb{T})$ with $||T_N^e\phi||_{H^1(\mathbb{R}^3 \times \mathbb{T})} \lesssim ||\phi||_{\dot{H}^1(\mathbb{R}^4)}$.
\begin{theorem}\label{theorem}
\noindent Assume that $\mathcal{O}=(N_k,t_k,x_k)_{k} \in \tilde{\mathcal{F}}_e$, $\phi \in \dot{H^1}(\mathbb{R}^4)$, and let $U_k(0)=\Pi_{t_k,x_k}(T^e_{N_k}\phi)$: \vspace{3mm}

\noindent (1) For $k$ large enough, there is a nonlinear solution $U_k \in X^1(\mathbb{R})$ of the equation (1.2) satisfying: 
\begin{equation}\label{euqation}
||U_k||_{X^1(\mathbb{R})} \lesssim_{E_{\mathbb{R}^4}(\phi)} 1.
\end{equation}
\noindent (2) There exists a Euclidean solution $u\in C(\mathbb{R}:\dot{H}^1(\mathbb{R}^4))$ of 
\begin{equation}\label{euqation} (i\partial_t+\Delta_{\mathbb{R}^4})u=|u|^2u
\end{equation}
\noindent with scattering data $\phi^{\pm  \infty}$ defined as in (5.4) such that up to a subsequence: for any $\epsilon > 0$, there exists $T(\phi,\epsilon)$ such that for all $T \geq T(\phi,\epsilon)$ there exists $R(\phi,\epsilon,T)$ such that for all $R \geq R(\phi,\epsilon,T)$, there holds that 
\begin{equation}\label{euqation} ||U_k-\tilde{u}_k||_{X^1(\{ |t-t_k| \leq TN_k^{-2}\})} \leq \epsilon,
\end{equation}
\noindent for $k$ large enough, where 
\begin{equation}\label{equation}
(\pi_{-x_k}\tilde{u})(x,t)=N_k \eta(N_k \Psi^{-1}(x)/R)u(N_k \Psi^{-1}(x),N_k^2(t-t_k)).
\end{equation}
\noindent In addition, up to a subsequence, 
\begin{equation} ||U_k(t)-\Pi_{(t_k-t,x_k)}T^e_{N_k} \phi^{\pm  \infty}||_{X^1(\{ |t-t_k| \leq TN_k^{-2}\})} \leq \epsilon
\end{equation}
\noindent for $k$ large enough (depending on ($\phi,\epsilon,T,R$)).
\end{theorem}
\noindent \emph{Proof:} This theorem is the analogue of [38, Theorem 5.4] (see also [14, 20]) and the proofs are similar so we omit it.\vspace{3mm}

\noindent \textbf{For quintic $\mathbb{R}^2\times \mathbb{T}$ problem:} \vspace{3mm}

\noindent We can also consider Euclidean profiles for the 3 dimensional case as follows. Since Theorem 5.6, Lemma 5.7 and Theorem 5.8 are analogues of Theorem 5.2, Lemma 5.3 and Theorem 5.4 respectively and the proofs are similar, we will omit the proofs except for the scattering norm control in Lemma 5.7. (see also [14, section 5], [20, section 5] and [38, section 5])\vspace{3mm}

\noindent \textbf{Euclidean Profiles.} The Euclidean profiles define a regime where we can compare solutions of cubic NLS on $\mathbb{R}^3$ with those on $\mathbb{R}^2\times \mathbb{T}$. We fix a spherically symmetric function $\eta \in C_{0}^{\infty}(\mathbb{R}^3)$ supported in the ball of radius 2 and equal to 1 in the ball of radius 1. Given $\phi \in H^1(\mathbb{R}^3)$ and a real number $N \geq 1$, we define:
\[
Q_{N} \phi \in H^1(\mathbb{R}^3)   \quad  \quad  (Q_{N} \phi)  (x)=\eta(\frac{x}{N^{\frac{1}{2}}}) \phi(x)
\]
\begin{equation}\label{euqation}
\phi_{N}\in H^1(\mathbb{R}^3)  \quad  \quad \phi_{N}(x)=N^{\frac{1}{2}}(Q_N \phi)(Nx)
\end{equation}
\[
f_N\in H^1(\mathbb{R}^2 \times \mathbb{T}) \quad \quad f_N(y)=\phi_N(\Psi^{-1}(y))
\]
\noindent where $\Psi$ is the identity map from the unit ball of $\mathbb{R}^3$ to $\mathbb{R}^2\times \mathbb{T}$. Thus $ Q_N \phi $ is a compactly supported modification of the profile $\phi$, $\phi_{N}$ is an $\dot{H^{1}}$-invariant rescaling of $Q_N \phi $, and $f_{N}$ is the function obtained by transferring $\phi_{N}$ to a neighborhood of $0$ in $\mathbb{R}^2 \times \mathbb{T}$. Notice that 
\[
\|f_{N}\|_{H^1(\mathbb{R}^2 \times \mathbb{T})}  \lesssim \|\phi\|_{\dot{H^{1}}(\mathbb{R}^3)}.
\]
\noindent  and for this case, we can apply the scattering result for 3d energy critical NLS by J. Colliander, M. Keel, G. Staffilani, H. Takaoka and T. Tao ([9]) in the following form:
\begin{theorem}\label{theorem}
Assume $\psi \in \dot{H^{1}}(\mathbb{R}^3)$, then there is a unique global solution $v \in C(\mathbb{R}:\dot{H^{1}}(\mathbb{R}^3))$ of the initial-value problem 
\begin{equation}\label{euqation}
(i\partial_{t} +\Delta_{\mathbb{R}^3})v=v|v|^{4}, \quad v(0)=\psi,
\end{equation}
\noindent and
\begin{equation}\label{euqation}
\||\nabla_{\mathbb{R}^3}v|\|_{(L_{t}^{\infty}L_{x}^{2} \cap L_{t}^{2}L_{x}^{6})(\mathbb{R}^3\times \mathbb{R})} \leq \tilde{C}(E_{\mathbb{R}^3}(\psi)).
\end{equation}
\noindent Moreover this solution scatters in the sense that there exists $\psi^{\pm \infty} \in \dot{H^{1}}(\mathbb{R}^3)$ such that 
\begin{equation}\label{euqation}
\|v(t)-e^{it\Delta} \psi^{\pm \infty}\|_{\dot{H^{1}}(\mathbb{R}^3)} \to 0
\end{equation}
as $t \to \pm \infty$. Besides if $\psi \in H^{5}(\mathbb{R}^3)$, then $v \in C(\mathbb{R}:H^5(\mathbb{R}^3))$ and 
\[
\sup\limits_{t\in \mathbb{R}} \|v(t)\|_{H^{5}(\mathbb{R}^3)}  \lesssim_{\|\psi\|_{H^{5}(\mathbb{R}^3)}}   1.
\] \end{theorem}
\noindent Based on the above result, we have:
\begin{theorem}\label{theorem}
\noindent Assume $\phi \in \dot{H^{1}}(\mathbb{R}^3)$, $T_{0} \in (0,\infty)$, and $\rho \in \{0,1\}$ are given, and we define $f_{N}$ as before. Then the following conclusions hold:\vspace{3mm} 

\noindent $(1)$ There is $N_0=N_0 (\phi,T_0)$ sufficiently large such that for any $N\geq N_0$, there is a unique solution $U_{N} \in C((-T_0 N^{-2},T_0 N^{-2}); H^1(\mathbb{R}^2\times \mathbb{T}))$ of the initial-value problem 
\begin{equation}\label{euqation}
(i\partial_t+\Delta)U_{N}=\rho U_{N} |U_{N}|^{4} ,\quad and \quad  U_{N}(0)=f_{N}.
\end{equation}
\noindent Moreover, for any $N \geq N_{0}$, 
\begin{equation}\label{euqation}
\|U_N\|_{X^{1}(-T_0 N^{-2},T_0 N^{-2})} \lesssim_{E_{\mathbb{R}^3}(\phi)} 1.
\end{equation}
$(2)$ Assume $\epsilon_{1} \in (0,1]$ is sufficiently small (depending only on $E_{\mathbb{R}^3}(\phi)$), $\phi ^{'} \in H^{5}(\mathbb{R}^3)$, and $\|\phi-\phi^{'}\|_{\dot{H^{1}}(\mathbb{R}^3)} \leq \epsilon_{1}$. Let $v^{'} \in C(\mathbb{R}:H^5)$ denote the solution of the initial-value problem 
\[
(i\partial_t+\Delta_{\mathbb{R}^3}) v^{'}=\rho v^{'} |v^{'}|^{4}, \quad v^{'} (0)=\phi^{'}.
\] 
\noindent For $R \geq 1$ and $N \geq 10R$, we define 
\[ v_{R}^{'}(x,t)=\eta(\frac{x}{R})v^{'}(x,t)    \quad \quad (x,t) \in \mathbb{R}^{3}\times (-T_0,T_0),
\]
\begin{equation}\label{euqation}
v_{R,N}^{'}(x,t)=N^{\frac{1}{2}} v_{R}^{'}(Nx,N^2 t) \quad \quad (x,t) \in \mathbb{R}^{3}\times (-T_0 N^{-2},T_0 N^{-2}),
\end{equation}
\[
V_{R,N}(y,t)=v_{R,N}^{'}(\Psi^{-1}(y),t) \quad \quad (y,t) \in \mathbb{R}^2\times \mathbb{T} \times (-T_0 N^{-2},T_0 N^{-2}).
\]
\noindent Then there is $R_0 \geq 1$ (depending on $T_0$ and $\phi^{'}$ and $\epsilon_1$) such that, for any $R \geq R_{0}$ and $N \geq 10R$,
\begin{equation}\label{euqation}
\limsup\limits_{N \to \infty} \|U_N-V_{R,N}\|_{X^{1}(-T_0 N^{-2},T_0 N^{-2})}  \lesssim_{E_{\mathbb{R}^3}(\phi)} \epsilon_1.
\end{equation} 
\end{theorem}

\begin{lemma}Suppose that $\phi \in \dot{H^1}(\mathbb{R}^3), \epsilon > 0$, and  $I \subset \mathbb{R} $ is an interval. Assume that
\begin{equation}\label{euqation}
||\phi||_{\dot{H^1} (\mathbb{R}^3)} \leq 1, \quad ||\nabla_{x} e^{it\Delta}\phi||_{L_t^2 L_x^6(\mathbb{R}^3 \times I)} \leq \epsilon.
\end{equation}
\noindent For $N \geq 1$, we define as before: 
\[Q_N \phi=\eta(N^{-1/2}x)\phi(x),\quad \phi_N=N^{\frac{1}{2}}(Q_N\phi)(Nx),\quad f_N(y)=\phi_N(\Psi^{-1}(y)).
\]
\noindent Then there exists $N_0=N_0(\phi,\epsilon)$ such that for any $N \geq N_0$,
\[ ||e^{it\Delta} f_{N}||_{Z(N^{-2}I)} \lesssim \epsilon.
\]
\end{lemma}
\noindent \emph{Proof.} By using similar analysis as in [14, Lemma 5.3] (see also [38, Lemma 5.3] and Lemma 5.3 in this paper), it suffices to prove: there exists $T_0$ such that for any $N>1$: 
\begin{equation}\label{euqation}
||e^{it\Delta} f_N||_{Z(\mathbb{R} \setminus (-N^{-2}T_0,N^{-2}T_0))} \lesssim \epsilon
\end{equation}
\noindent and we can obtain:
\begin{equation}\label{euqation} 
(\sum_{M} M^{\frac{1}{4}}||e^{it\Delta}P_M f_N||^2_{L_{x,y,t}^{8}(\mathbb{R}^2\times \mathbb{T}\times \{N^2|t|\geq T\})})^{\frac{1}{2}}\lesssim T^{-\frac{5}{8}} .
\end{equation}
\noindent Also, by using Stricharz estimate, for $p>10$, we have
\begin{equation}
(\sum_{M} M^{\frac{10}{p}-1}||e^{it\Delta}P_M f_N||^2_{l_{\gamma}^{\frac{4p}{p-4}}L_{x,y,t}^{p}(\mathbb{R}^2\times \mathbb{T} \times \mathbb{R})})^{\frac{1}{2}} \lesssim ||\phi_{N}||_{H^1} \lesssim 1.
\end{equation}
\noindent By interpolation, we can obtain (5.33). That finishes the proof of Lemma 5.7.\vspace{5mm}
\begin{theorem}\label{theorem}
\noindent Assume that $\mathcal{O}=(N_k,t_k,x_k)_{k} \in \tilde{\mathcal{F}}_e$, $\phi \in \dot{H^1}(\mathbb{R}^3)$, and let $U_k(0)=\Pi_{t_k,x_k}(T^e_{N_k}\phi)$: \vspace{3mm}

\noindent (1) For $k$ large enough, there is a nonlinear solution $U_k \in X^1(\mathbb{R})$ of the equation (1.3) satisfying: 
\begin{equation}\label{euqation}
||U_k||_{X^1(\mathbb{R})} \lesssim_{E_{\mathbb{R}^3}(\phi)} 1.
\end{equation}
\noindent (2) There exists a Euclidean solution $u\in C(\mathbb{R}:\dot{H}^1(\mathbb{R}^3))$ of 
\begin{equation}\label{euqation} (i\partial_t+\Delta_{\mathbb{R}^3})u=|u|^4u
\end{equation}
\noindent with scattering data $\phi^{\pm  \infty}$ defined as in (5.4) such that up to a subsequence: for any $\epsilon > 0$, there exists $T(\phi,\epsilon)$ such that for all $T \geq T(\phi,\epsilon)$ there exists $R(\phi,\epsilon,T)$ such that for all $R \geq R(\phi,\epsilon,T)$, there holds that 
\begin{equation}\label{euqation} ||U_k-\tilde{u}_k||_{X^1(\{ |t-t_k| \leq TN_k^{-2}\})} \leq \epsilon,
\end{equation}
\noindent for $k$ large enough, where 
\begin{equation}\label{equation}
(\pi_{-x_k}\tilde{u})(x,t)=N_k^{\frac{1}{2}} \eta(N_k \Psi^{-1}(x)/R)u(N_k \Psi^{-1}(x),N_k^2(t-t_k)).
\end{equation}
\noindent In addition, up to a subsequence, 
\begin{equation} ||U_k(t)-\Pi_{(t_k-t,x_k)}T^e_{N_k} \phi^{\pm  \infty}||_{X^1(\{ |t-t_k| \leq TN_k^{-2}\})} \leq \epsilon
\end{equation}
\noindent for $k$ large enough (depending on ($\phi,\epsilon,T,R$)).
\end{theorem}

\section{Profile Decomposition}
\begin{definition} [Frames and Profiles](1) We define a frame to be sequence $(N_k,t_k,p_k)_{k} \in 2^{\mathbb{Z}} \times \mathbb{R} \times (\mathbb{R}^3 \times \mathbb{T}) $ . And we can define some types of profiles as follows.

\noindent a) A Euclidean frame is a sequence $\mathcal{F}_e =(N_k,t_k,p_k)$ with $N_k \geq 1, N_k \rightarrow \infty, t_k \in \mathbb{R} , p_k \in \mathbb{R}^3 \times \mathbb{T}$.

\noindent b) A  Scale-one frame is a sequence $\mathcal{F}_1 =(1,t_k,p_k)$ with $t_k \in \mathbb{R}, p_k \in \mathbb{R}^3 \times \mathbb{T}$.

\noindent (2) We say that two frames $(N_k,t_k,p_k)_{k}$ and $(M_k,s_k,q_k)_{k}$ are orthogonal if
\[ \lim_{k \rightarrow +\infty} (|ln \frac{N_k}{M_k}|+N_k^2|t_k-s_k|+N_k|(p_k-q_k)|)=+\infty.
\]
\noindent(3) We associate a profile defined as: 

\noindent a) If $\mathcal{O}=(N_k,t_k,p_k)_k $ is a Euclidean frame and for $\phi \in \dot{H}^1(\mathbb{R}^4)$ we define the Euclidean profile associated to $(\phi,\mathcal{O})$ as the sequence $\tilde{\phi}_{\mathcal{O},k}$ with
\[\tilde{\phi}_{\mathcal{O},k}=\Pi_{t_k,p_k}(T^{e}_{N_k}\phi)(x,y).
\]

\noindent b) If $\mathcal{O}=(1,t_k,p_k)_k$  is a scale one frame, if $W \in H^1(\mathbb{R}^3 \times \mathbb{T})$, we define the scale one profile associated to $(W,\mathcal{O})$ as $\tilde{W}_{\mathcal{O},k}$ with
\[\tilde{W}_{\mathcal{O},k}=\Pi_{t_k,p_k} W.
\]

\noindent (4) Finally, we say that a sequence of functions $\{f_k\}_k \subset H^1(\mathbb{R}^3 \times \mathbb{T})$ is absent from a frame $\mathcal{O}$ if, up to a subsequence: 

\noindent $\langle f_k,\tilde{\psi}_{\mathcal{O},k} \rangle_{H^1 \times H^1}\rightarrow 0$ as $k \rightarrow \infty$ for any profile $\tilde{\psi}_{\mathcal{O},k}$ associated with $\mathcal{O}$.
\end{definition}
\noindent \emph{Remark 1.} The definition for the case of quintic $\mathbb{R}^2\times \mathbb{T}$ is quite similar. We can replace $\mathbb{R}^3\times \mathbb{T}$ by $\mathbb{R}^2\times \mathbb{T}$ in the Definition 6.1 so we will not repeat it again.\vspace{3mm}

\noindent \emph{Remark 2.} It is very convenient to use the language of frames and profiles to unify Euclidean profiles and scale-one profiles. There are some useful properties about the equivalence of frames. We refer [14, 18, 20, 38] for more information.\vspace{3mm}

\noindent \textbf{For cubic $\mathbb{R}^3\times \mathbb{T}$ problem:} \vspace{3mm}

\noindent Now, let us state a core lemma which is a key step of the main theorem in this section (profile decomposition). First, for a bounded sequence of functions $\{f_k\}$ in $H^1(\mathbb{R}^3 \times \mathbb{T})$, we define the following functional based on a Besov norm:
\begin{equation}\label{equation} 
\Lambda_{\infty}(\{f_k\})=\limsup\limits_{k \rightarrow \infty}||e^{it\Delta}f_k||_{L_t^{\infty}B_{\infty,\infty}^{-1}}=\limsup\limits_{k \rightarrow \infty} \sup\limits_{\{N,t,x,y\}} N^{-1} |(e^{it\Delta}P_N f_k)(x,y)| .
\end{equation}
\noindent Given a uniformly bounded sequence in $H^1$, we will extract some profiles whose Besov norms are big to ensure the Besov norm of the linear propagation of the remainder flow is small.
\begin{lemma}Let $v > 0$. Assume that $\phi_k$ is a sequence satisfying $||\phi_k||_{H^1(\mathbb{R}^3\times \mathbb{T})}<E$, then there exists a subsequence of $\phi_k$ (for convenience, we still use $\phi_k$), $A$ Euclidean profiles $\tilde{\varphi^{\alpha}}_{\mathcal{O}^\alpha,k}$, and $A$ scale-one profiles $\tilde{W^{\beta}}_{\mathcal{O}^\beta,k}$ such that, for any $k \geq 0$ in the subsequence
\begin{equation}\label{equation} 
\phi_k^{'}(x,y)=\phi_k(x,y)-\sum_{1\leq \alpha \leq A} \tilde{\varphi}^{\alpha}_{\mathcal{O}^\alpha,k}-\sum_{1\leq \beta \leq A}\tilde{W}^{\beta}_{\mathcal{O}^\beta,k} \end{equation} \noindent satisfies
\begin{equation}\label{equation} 
\Lambda_{\infty}(\{ \phi_k^{'}\}) < v.
\end{equation}
\noindent Besides, all the frames involved are pairwise orthogonal and $\phi_k^{'}$ is absent from all these frames.
\end{lemma}
\noindent \emph{Proof.} The proof is similar to [14, Lemma 6.5], [20, Lemma 5.4] and [38, Lemma 6.4]. We omit it.

\begin{theorem}[Profile decomposition]
Assume $\{ \phi_k\}_{k}$ is a sequence of functions satisfying
$||\phi_k||_{H^1(\mathbb{R}^3\times \mathbb{T})}<E$, up to a subsequence, then there exists a sequence of Euclidean profiles $\tilde{\varphi^{\alpha}}_{\mathcal{O}^\alpha,k}$, and scale-one profiles $\tilde{W^{\beta}}_{\mathcal{O}^\beta,k}$ such that, for any $J \geq 0$ 
\begin{equation}\label{equation} 
\phi_k(x,y)=\sum_{1\leq \alpha \leq J} \tilde{\varphi}^{\alpha}_{\mathcal{O}^\alpha,k}+\sum_{1\leq \beta \leq J}\tilde{W}^{\beta}_{\mathcal{O}^\beta,k}+R^J_k \end{equation} 
\noindent where $R^J_k$ is absent from the frames $\mathcal{O}^{\alpha}$ and satisfies
\begin{equation}
\limsup_{J\rightarrow \infty} \Lambda_{\infty}(\{R^J_k\})=0
\end{equation}
\noindent Additionally, we have the following orthogonal relation
\[||\phi_k||^2_{L^2}=\sum_{\alpha} || \tilde{\varphi^{\alpha}}_{\mathcal{O}^\alpha,k}||^2_{L^2}+\sum_{\beta} || \tilde{W}^{\beta}_{\mathcal{O}^\beta,k}||^2_{L^2}+||R_k^J||_{L^2}^2+o_k(1),
\]
\begin{equation}
 ||\nabla \phi_k||^2_{L^2}=\sum_{\alpha} ||\nabla \tilde{\varphi}^{\alpha}_{\mathcal{O}^\alpha,k}||^2_{L^2}+ \sum_{\beta} || \nabla \tilde{W}^{\beta}_{\mathcal{O}^\beta,k}||^2_{L^2}+||\nabla R_k^J||_{L^2}^2+o_k(1),
\end{equation}
\[|| \phi_k||^4_{L^4}=\sum_{\alpha} || \tilde{\varphi}^{\alpha}_{\mathcal{O}^\alpha,k}||^4_{L^4}+\sum_{\beta} ||\tilde{W}^{\beta}_{\mathcal{O}^\beta,k}||^4_{L^4}+o_{J,k}(1),
\]
\end{theorem}
\noindent \emph{Proof:} The proof is similar as in [20, Proposition 5.5]. We omit it. It mainly follows from Lemma 6.1, Lemma 2.2 and some frame equivalence properties. \vspace{3mm}

\noindent The following lemma is crucial, which explains how to use the Bosov norm to control the scattering norm of linear Schr{\"o}dinger propagation.
\begin{lemma} Assume a sequence $\{f_k\}_k$ satisfying $\sup_k ||f_k||_{H^1}<E$, we have the following estimate:
\begin{equation}
\limsup\limits_{k \rightarrow \infty} ||e^{it\Delta}f_k||_{Z(\mathbb{R})}\lesssim_{E} (\Lambda_{\infty}(\{f_k\}))^{\delta},
\end{equation}
\noindent where $\delta$ is some positive constant.
\end{lemma}
\noindent \emph{Proof.} By using interpolation and Stricharz estimate, 
\begin{align*}
||e^{it\Delta}f_k||^{p_0}_{Z(\mathbb{R})}&=\sum_{N\geq 1}N^{6-p_0} \|1_I(t) P_N u\|_{l^{\frac{2p_0}{p_0-3}}L_{x,y,t}^{p_0}(\mathbb{R}^3\times \mathbb{T} \times I_{\gamma})}^{p_0} \\
&\lesssim \sum_N (N^{-1}||P_N e^{it\Delta}f_k||_{L_{x,t}^{\infty}})^{\frac{p_0}{3}}(N^{\frac{9}{p_0}-1}||P_N e^{it\Delta}f_k||_{l^{\frac{4p_0}{3p_0-9}}L^{\frac{2p_0}{3}}})^{\frac{2p_0}{3}} \\
&\lesssim (\sup\limits_N N^{-1}||P_N e^{it\Delta}f_k||_{L_{x,t}^{\infty}})^{\frac{p_0}{3}} \sum_N (N||f_k||_{L^2})^{\frac{2p_0}{3}} \\
&\lesssim (\sup\limits_N N^{-1}||P_N e^{it\Delta}f_k||_{L_{x,t}^{\infty}})^{\frac{p_0}{3}} ||f_k||_{H^1}^{\frac{2p_0}{3}} \\
&\lesssim_E (\sup\limits_N N^{-1}||P_N e^{it\Delta}f_k||_{L_{x,t}^{\infty}})^{\frac{p_0}{3}}.
\end{align*}
\noindent That finishes the proof of Lemma 6.4 noticing that $||f_k||_{H^1}$ is uniformly bounded by $E$.\vspace{3mm}

\noindent \emph{Remark 1.} The range of $p_0$ ensures $\frac{22}{7}<\frac{2p_0}{3}<6$ so that we can use Stricharz estimate (Theorem 3.1) in the proof.\vspace{3mm}

\noindent \emph{Remark 2.} If we replace $f_k$ by $R_k^J$ which is the remainder term in the profile decomposition, we can obtain the control of the scattering norm of linear Schr{\"o}dinger propagation of the remainder flow.\vspace{3mm}

\noindent \emph{Remark 3.} We are mainly inspired by [14] and [20]. In those papers, the authors have used similar techniques to estimate the remainder flow.\vspace{3mm}

\noindent \textbf{For quintic $\mathbb{R}^2\times \mathbb{T}$ problem:} \vspace{3mm}

\noindent First, we state a core lemma which is a key step of the main theorem in this section (profile decomposition). First, for a bounded sequence of functions $\{f_k\}$ in $H^1(\mathbb{R}^2 \times \mathbb{T})$, we define the following functional based on a Besov norm:
\begin{equation}\label{equation} 
\Lambda_{\infty}(\{f_k\})=\limsup\limits_{k \rightarrow \infty}||e^{it\Delta}f_k||_{L_t^{\infty}B_{\infty,\infty}^{-\frac{1}{2}}}=\limsup\limits_{k \rightarrow \infty} \sup\limits_{\{N,t,x,y\}} N^{-\frac{1}{2}} |(e^{it\Delta}P_N f_k)(x,y)| .
\end{equation}
\noindent Given a uniformly bounded sequence in $H^1$, we will extract some profiles whose Besov norms are big to ensure the Besov norm of the linear propagation of the remainder flow is small.
\begin{lemma}Let $v > 0$. Assume that $\phi_k$ is a sequence satisfying $||\phi_k||_{H^1(\mathbb{R}^2\times \mathbb{T})}<E$, then there exists a subsequence of $\phi_k$ (for convenience, we still use $\phi_k$), $A$ Euclidean profiles $\tilde{\varphi^{\alpha}}_{\mathcal{O}^\alpha,k}$, and $A$ scale-one profiles $\tilde{W^{\beta}}_{\mathcal{O}^\beta,k}$ such that, for any $k \geq 0$ in the subsequence
\begin{equation}\label{equation} 
\phi_k^{'}(x,y)=\phi_k(x,y)-\sum_{1\leq \alpha \leq A} \tilde{\varphi}^{\alpha}_{\mathcal{O}^\alpha,k}-\sum_{1\leq \beta \leq A}\tilde{W}^{\beta}_{\mathcal{O}^\beta,k} \end{equation} \noindent satisfies
\begin{equation}\label{equation} 
\Lambda_{\infty}(\{ \phi_k^{'}\}) < v.
\end{equation}
\noindent Besides, all the frames involved are pairwise orthogonal and $\phi_k^{'}$ is absent from all these frames.
\end{lemma}

\begin{theorem}[Profile decomposition]
Assume $\{ \phi_k\}_{k}$ is a sequence of functions satisfying
$||\phi_k||_{H^1(\mathbb{R}^2\times \mathbb{T})}<E$, up to a subsequence, then there exists a sequence of Euclidean profiles $\tilde{\varphi^{\alpha}}_{\mathcal{O}^\alpha,k}$, and scale-one profiles $\tilde{W^{\beta}}_{\mathcal{O}^\beta,k}$ such that, for any $J \geq 0$ 
\begin{equation}\label{equation} 
\phi_k(x,y)=\sum_{1\leq \alpha \leq J} \tilde{\varphi}^{\alpha}_{\mathcal{O}^\alpha,k}+\sum_{1\leq \beta \leq J}\tilde{W}^{\beta}_{\mathcal{O}^\beta,k}+R^J_k \end{equation} 
\noindent where $R^J_k$ is absent from the frames $\mathcal{O}^{\alpha}$ and satisfies
\begin{equation}
\limsup_{J\rightarrow \infty} \Lambda_{\infty}(\{R^J_k\})=0.
\end{equation}
\noindent Additionally, we have the following orthogonal relation
\[||\phi_k||^2_{L^2}=\sum_{\alpha} || \tilde{\varphi^{\alpha}}_{\mathcal{O}^\alpha,k}||^2_{L^2}+\sum_{\beta} || \tilde{W}^{\beta}_{\mathcal{O}^\beta,k}||^2_{L^2}+||R_k^J||_{L^2}^2+o_k(1),
\]
\begin{equation}
 ||\nabla \phi_k||^2_{L^2}=\sum_{\alpha} ||\nabla \tilde{\varphi}^{\alpha}_{\mathcal{O}^\alpha,k}||^2_{L^2}+ \sum_{\beta} || \nabla \tilde{W}^{\beta}_{\mathcal{O}^\beta,k}||^2_{L^2}+||\nabla R_k^J||_{L^2}^2+o_k(1),
\end{equation}
\[|| \phi_k||^6_{L^6}=\sum_{\alpha} || \tilde{\varphi}^{\alpha}_{\mathcal{O}^\alpha,k}||^6_{L^6}+\sum_{\beta} ||\tilde{W}^{\beta}_{\mathcal{O}^\beta,k}||^6_{L^6}+o_{J,k}(1).
\]
\end{theorem}
\noindent The following estimate is crucial and it is the analogue of Lemma 6.3.

\begin{lemma}[Control of the remainder term]Assume a sequence $\{f_k\}_k$ satisfying $\sup_k ||f_k||_{H^1}<E$, we have the following estimate:
\begin{equation}
\limsup\limits_{k \rightarrow \infty} ||e^{it\Delta}f_k||_{Z(\mathbb{R})}\lesssim_{E} (\Lambda_{\infty}(\{f_k\}))^{\delta},
\end{equation}
\noindent where $\delta$ is some positive constant.
\end{lemma}
\noindent \emph{Proof.} By using interpolation and Stricharz estimate, 
\begin{align*}
||e^{it\Delta}f_k||_{Z(\mathbb{R})}&=\sum_{p_0=8,10}(\sum_{N\geq 1}N^{5-\frac{p_0}{2} } \|1_I(t) P_N u\|_{l^{\frac{4p_0}{p_0-4}}L_{x,y,t}^{p_0}(\mathbb{R}^2\times \mathbb{T} \times I_{\gamma})}^{p_0})^\frac{1}{p_0} \\
&\lesssim \sum_{p_0=8,10} \sum_N (N^{-\frac{1}{2}}||P_N e^{it\Delta}f_k||_{L_{x,t}^{\infty}})^{\frac{1}{2}}(N^{\frac{10}{p_0}-\frac{1}{2}}||P_N e^{it\Delta}f_k||_{l^{\frac{2p_0}{p_0-4}}L^{\frac{p_0}{2}}})^{\frac{1}{2}} \\
&\lesssim (\sup\limits_N N^{-\frac{1}{2}}||P_N e^{it\Delta}f_k||_{L_{x,t}^{\infty}})^{\frac{1}{2}} \sum_N (N||f_k||_{L^2})^{\frac{1}{2}} \\
&\lesssim (\sup\limits_N N^{-\frac{1}{2}}||P_N e^{it\Delta}f_k||_{L_{x,t}^{\infty}})^{\frac{1}{2}} ||f_k||_{H^1}^{\frac{1}{2}} \\
&\lesssim_E (\sup\limits_N N^{-\frac{1}{2}}||P_N e^{it\Delta}f_k||_{L_{x,t}^{\infty}})^{\frac{1}{2}}.
\end{align*}
\noindent That finishes the proof of Lemma 6.7 noticing that $||f_k||_{H^1}$ is uniformly bounded by $E$.

\section{Induction on Energy}
\noindent We are now ready to prove the main theorem. We follow an induction on energy method formalized in [21, 22]. Similar as related results ([7, 14, 20, 38]), we also consider the full energy since it is a $H^1$ data problem. We define the following functional
\[ \Lambda(L)=\textmd{sup}\{||u||_{Z(I)}:u \in X^1_{loc}(I),E(u)+M(u) \leq L\}
\]
\noindent where the supremum is taken over all strong solutions of full energy less than $L$. According to the local theory, this is sublinear in $L$ and finite for $L$ sufficiently small. We also define
\[L_{max}= \textmd{sup}\{ L: \Lambda(L) < +\infty \}.
\]\noindent In other to prove the large data scattering of (1.2) and (1.3), it suffices to show that $ L_{max}=+\infty$ according to Theorem 4.4 and Theorem 4.9. That is our goal. The key proposition is as follows (Theorem 7.1 and Theorem 7.4):\vspace{3mm}

\noindent \textbf{For cubic $\mathbb{R}^3\times \mathbb{T}$ problem:} 

\begin{theorem}\label{theorem} Assume that $L_{max} < +\infty $. Let $\{ t_{k}\}_{k} ,\{ a_{k}\}_{k} ,\{ b_{k}\}_{k} $ be arbitrary sequences of real numbers and $\{u_k\}$ be a sequence of solutions to (1.1) such that $u_{k} \in X^{1}_{c,loc}(t_k-a_k,t_k+b_k)$ and satisfying
\begin{equation}\label{equation} L(u_k) \rightarrow L_{max},\quad ||u_k||_{Z(t_k-a_k,t_k)} \rightarrow +\infty,\quad ||u_k||_{Z(t_k,t_k+b_k)} \rightarrow +\infty. \end{equation}
Then passing to a subsequence, there exists a sequence $x_{k} \in \mathbb{R}^3 $ and $\omega \in H^{1}(\mathbb{R}^3 \times \mathbb{T})$ such that
\begin{equation}\label{equation} \omega_k(x,y) =u_k(x-x_k,y,t_k) \rightarrow \omega \end{equation}
strongly in $H^{1}(\mathbb{R}^3 \times \mathbb{T})$.
\end{theorem}
\noindent The proof mainly follows from profile decomposition and perturbation theory. We will give the proof at the end of Section 7. Now let us show how to use Theorem 7.1 to close the contradiction argument and finish the proof of the main theorem. The following analysis (Corollary 7.2 and Theorem 7.3) is similar to [14] (see also [38]).
\begin{corollary} Assume that $L_{max} < +\infty $. Then there exists $u \in X^1_{c,loc}(\mathbb{R})$ solving (1.1) and a Lipschitz function $\underline{x}:\mathbb{R} \rightarrow \mathbb{R}^3$ such that $L(u)=L_{max}$ and
\begin{equation}\label{equation} \sup\limits_{t \in \mathbb{R}} |\underline{x}^{'}(t)| \lesssim 1, \end{equation} 
\[ (u(x-\underline{x}(t),y,t) :t \in \mathbb{R})\quad \textmd{is precompact in} \quad H^1(\mathbb{R}^3 \times \mathbb{T}).
\]
\end{corollary}
\noindent \emph{Remark.} The proof is similar to [14, Corollary 7.2] by using perturbation theory. We omit it.
\begin{theorem}\label{theorem}Assume that $u$ satisfies the conclusion of Corollary 7.2, then $u=0$. 
\end{theorem}
\noindent \emph{Proof:} Assume $u \neq 0$. Then, from the compactness property, we see that there exists $\rho >0$ such that 
\begin{equation}\label{equation} \inf\limits_{t \in \mathbb{R}}min (||u(t)||_{L^4_{x,y}(\mathbb{R}^3 \times \mathbb{T})},||u(t)||_{L^2_{x,y}(\mathbb{R}^3 \times \mathbb{T})}) \geq \rho.
\end{equation}
\noindent Now let us consider the conserved momentum 
\[ P(u)=Im \int_{\mathbb{R}^3\times \mathbb{T}}\bar{u}(x,y,t) \partial_{x_1} u(x,y,t) dxdy.
\]
\noindent Considering the Galilean transform 
\[v(z,t)=e^{-i|\xi_0|^2t+i \langle z,\xi_0 \rangle} u(z-2\xi_0 t,t),
\]
\noindent and letting
\[\xi_0=-\frac{P(u)}{M(u)},
\]
\noindent without loss of generality, we can assume that 
\begin{equation}\label{equation} P(u)=0.
\end{equation} \noindent Then we define the Virial action by 
\[ A_{R}(t)=\int_{\mathbb{R}^3\times \mathbb{T}}  \chi_{R}(x_1-\underline{x}_1(t))(x_1-\underline{x}_1(t))Im[\bar{u}(x,y,t)\partial_{x_1} u(x,y,t)]dxdy
\]
\noindent for $\chi_R(x)=\chi(x/R)$ and $\chi$ satisfies $\chi(x)=1$ when $|x| \leq 1$ and $\chi(x)=0$ when $|x| \geq 2$ ($x \in \mathbb{R}$).\vspace{3mm}

\noindent On one hand, clearly 
\begin{equation}\label{equation} \sup\limits_{t}|A_{R}(t)| \lesssim R .
\end{equation}
\noindent On the other hand, we compute that
\begin{align*}
\frac{d}{dt}A_R &=- \underline{x}_1^{'}(t) Im \int_{\mathbb{R}^3\times \mathbb{T}}\bar{u}(x,y,t) \partial_{x_1}  u(x,y,t) dxdy \\
& - \underline{x}_1^{'}(t)\int_{\mathbb{R}^3\times \mathbb{T}} \{(\chi ^{'})_R(x_1-\underline{x}_1(t))\frac{x_1-\underline{x}_1(t)}{R} -(1- \chi_{R}(x_1-\underline{x}_1(t)))\} Im [\bar{u}(x,y,t) \partial_{x_1}  u(x,y,t)] dxdy\\
&+  \int_{\mathbb{R}^3\times \mathbb{T}} \chi_{R}(x_1-\underline{x}_1(t))(x_1-\underline{x}_1(t)) \partial_t Im[\bar{u}(x,y,t)\partial_{x_1} u(x,y,t)]dxdy.
\end{align*}

\noindent The first term will vanish automatically based on the assumption (7.5) and the second term can be bounded by  
\[ \int_{\{|x-\underline{x}(t)| \geq R\}} \int_{\mathbb{T}}[|u(x,y,t)|^2+|\nabla u(x,y,t)|^2]dxdy=O_{R}(t),
\]
 $\qquad \qquad \qquad  \qquad  \qquad  \qquad  \sup\limits_t O_{R}(t) \rightarrow 0 \qquad$ as $\qquad  \quad R  \rightarrow +\infty $.

\noindent Notice that
\[  \partial_t Im[\bar{u}(x,y,t) \partial_{x_1} u(x,y,t)]=\partial_{x_1} \Delta \frac{|u|^2}{2}-2div\{Re[\partial_{x_1}\bar{u} \nabla u]\}-\frac{1}{4}\partial_{x_1} |u|^4.
\]
\noindent For the last term, we have 
\begin{align*}
\frac{d}{dt}A_R &=\int_{\mathbb{R}^3\times \mathbb{T}}\chi_{R}(x_1-\underline{x}_1(t))[\frac{1}{4}|u(x,y,t)|^4+\frac{1}{2}|\partial_{x_1} u(x,y,t)|^2]dxdy \\
&+ \int_{\mathbb{R}^3\times \mathbb{T}}\chi_{R}^{'}(x_1-\underline{x}_1(t)) \frac{x_1-\underline{x}_1(t)}{R}[\frac{1}{4}|u(x,y,t)|^4+\frac{1}{2}|\partial_{x_1} u(x,y,t)|^2]dxdy\\
&-  \int_{\mathbb{R}^3\times \mathbb{T}} \frac{|u(x,y,t)|^2}{2} \partial_{x_1}^3[ \chi_{R}(x_1-\underline{x}_1(t))(x_1-\underline{x}_1(t)) ]dxdy +O_R(t)\\
&= \int_{\mathbb{R}^3\times \mathbb{T}} [\frac{1}{4}|u(x,y,t)|^4+\frac{1}{2}|\partial_{x_1} u(x,y,t)|^2]dxdy+\tilde{O}_{R}(t).
\end{align*}
\noindent Integrating this equality, we obtain 
\[|A_R(t)-A_R(0)| \geq Ct\rho-t \sup\limits_t \tilde{O}_R(t).
\]
\noindent Taking $R$ sufficiently large enough, when $t$ is sufficiently large, we see there is a contradiction. This finishes the proof of Theorem 7.3.\vspace{3mm}

\noindent \textbf{Proof of Theorem 7.1:} The proof is similar to [20, Theorem 6.1] (see also [38, Theorem 7.1]). First, we apply profile decomposition to the bounded $H^1$ sequence $u_k(0)$ and then consider three cases, i.e. no profile, only one profile and multiple profiles. The first two cases can be easily handled by using Lemma 6.3 and approximation result and the last case can be handled by constructing approximate solution and perturbation theory.\vspace{3mm}

\noindent Without loss of generality, we assume $t_k=0$, and we apply Theorem 6.2 (Profile decomposition) to $\{u_k(0)\}_k$ which is a bounded sequence in $H^1(\mathbb{R}^3\times \mathbb{T})$. For all $J$, we have
\begin{equation}\label{equation} 
u_k(0)=\sum_{1\leq \alpha \leq J} \tilde{\varphi}^{\alpha}_{\mathcal{O}^\alpha,k}+\sum_{1\leq \beta \leq J}\tilde{W}^{\beta}_{\mathcal{O}^\beta,k}+R^J_k. \end{equation} 

\noindent \emph{Case 1:} There are no profiles. By using Lemma 6.3 (scattering norm estimate), we have, if we take $J$ sufficiently large, we will have:
\[ ||e^{it\Delta} u_{k}(0)||_{Z(\mathbb{R})} = ||e^{it\Delta} R^{J}_{k}||_{Z(\mathbb{R})}\leq  \delta_{0}/2
\]\noindent for k sufficiently large, where $\delta_{0}$ is given in Theorem 4.3. Then we know that $u_k$ can be extended on $\mathbb{R}$ and that
\[ \lim\limits_{k \rightarrow +\infty} ||u_k||_{Z(\mathbb{R})} \leq \delta_{0}.
\]

\noindent It is a contradiction. Hence there are at least one profile. There are two other cases left: only one profile (\emph{Case 2}) and multiple profiles (\emph{Case 3}). Furthermore, only one profile contains two cases, i.e. only one Euclidean profile (\emph{Case 2a}) and only scale-one profile (\emph{Case 2b}). Except for (\emph{Case 1}), we will also rule out \emph{Case 2a} and \emph{Case 3}. Actually the conclusion statement in Theorem 7.1 is the only case that will happen. That is our goal for Theorem 7.1 which is an important preparation for the compactness argument (Corollary 7.2). \vspace{3mm}

\noindent Moreover, for every linear profile, we define the associated nonlinear profile as the maximal solution of (1.2) with the corresponding initial data as in [14].\vspace{3mm}

\noindent For any profile, we consider operator $L$ such that
\[L(\alpha):=\lim\limits_{k\rightarrow \infty} (E(\tilde{\varphi}^{\alpha}_{\mathcal{O}^\alpha,k})+M(\tilde{\varphi}^{\alpha}_{\mathcal{O}^\alpha,k}))\in (0,L_{max}].
\]
\noindent According to the orthogonal properties in profile decomposition, we have:
\begin{equation}\label{equation} \lim\limits_{J \rightarrow +\infty}[\sum_{1 \leq \alpha,\beta \leq J}[L_{E}(\alpha)+L_{1}(\beta)]+\lim\limits_{k \rightarrow +\infty} L(R_k^J)] \leq L_{max}. \end{equation}
\noindent \emph{Case 2a:} There are only one Euclidean profile in the profile decomposition, that is
\[ u_k(0)=\tilde{\phi}_{\varepsilon,k}+o_k(1)
\]
\noindent in $H^1(\mathbb{R}^3 \times \mathbb{T})$, where $\varepsilon$ is a Euclidean frame. In this case, since the corresponding nonlinear profile $U_k$ satisfies $||U_k||_{Z(\mathbb{R})} \lesssim_{E_{\mathbb{R}^4}(\phi)} 1 $ and $  \lim\limits_{k \rightarrow +\infty } ||U_k(0)-u_k(0)||_{H^1(\mathbb{R}^3 \times \mathbb{T})} \rightarrow 0 $. We can use Theorem 4.5 to deduce that 
\[ ||u_k||_{Z(\mathbb{R})} \lesssim ||u_k||_{X^1(\mathbb{R})} \lesssim_{L_{max}} 1,
\]
\noindent which contradicts (7.1).\vspace{3mm}

\noindent \emph{Case 2b:} There are only one scale-one profile in the profile decomposition, we have that 
\[ u_k(0)=\tilde{\omega}_{\mathcal{O},k}+o_k(1)
\]
in $H^1(\mathbb{R}^3\times \mathbb{T})$, where $\mathcal{O}=\{1,t_k,x_k\}$ is a scale-one frame. If $t_k \equiv 0$, this is precisely the conclusion(7.2). \vspace{3mm}

\noindent If $t_k \rightarrow +\infty$, then 
\[||e^{it\Delta_{\mathbb{R}^3\times \mathbb{T}}}\tilde{\omega}_{\mathcal{O},k}||_{Z(a_k,0)} \leq ||e^{it\Delta_{\mathbb{R}^3\times \mathbb{T}}}\tilde{\omega}_{\mathcal{O},k}||_{Z(-\infty,0)}=||e^{it\Delta_{\mathbb{R}^3\times \mathbb{T}}}\omega||_{Z(-\infty,-t_k)}
\]
\noindent which goes to 0 as $t_k \rightarrow +\infty$. Using Theorem 4.3, we see that, for $k$ large enough, 
\[ ||u_k||_{Z(-\infty,0)} \leq \delta_0.
\]
\noindent It contradicts (7.1). The case $t_k \rightarrow  -\infty$ is similar.\vspace{3mm}

\noindent \emph{Case 3:} There are multiple profiles in the profile decomposition. We can construct approximate equation and use perturbation theory to rule out this case. It is similar as to [14, Proposition 7.1, case 3] and [20, Proposition 6.1, case 3] (see also [38, Theorem 7.1]). We omit it.\vspace{3mm}

\noindent \textbf{For quintic $\mathbb{R}^2\times \mathbb{T}$ problem:} we have, 

\begin{theorem}\label{theorem} Assume that $L_{max} < +\infty $. Let $\{ t_{k}\}_{k} ,\{ a_{k}\}_{k} ,\{ b_{k}\}_{k} $ be arbitrary sequences of real numbers and $\{u_k\}$ be a sequence of solutions to (1.1) such that $u_{k} \in X^{1}_{c,loc}(t_k-a_k,t_k+b_k)$ and satisfying
\begin{equation}\label{equation} L(u_k) \rightarrow L_{max},\quad ||u_k||_{Z(t_k-a_k,t_k)} \rightarrow +\infty,\quad ||u_k||_{Z(t_k,t_k+b_k)} \rightarrow +\infty. \end{equation}
Then passing to a subsequence, there exists a sequence $x_{k} \in \mathbb{R}^2 $ and $\omega \in H^{1}(\mathbb{R}^2 \times \mathbb{T})$ such that
\begin{equation}\label{equation} \omega_k(x,y) =u_k(x-x_k,y,t_k) \rightarrow \omega \end{equation}
strongly in $H^{1}(\mathbb{R}^2 \times \mathbb{T})$.
\end{theorem}
\noindent The proof of Theorem 7.4 mainly follows from profile decomposition and perturbation theory. Based on this, we can prove:
\begin{corollary} Assume that $L_{max} < +\infty $. Then there exists $u \in X^1_{c,loc}(\mathbb{R})$ solving (1.1) and a Lipschitz function $\underline{x}:\mathbb{R} \rightarrow \mathbb{R}^2$ such that $L(u)=L_{max}$ and
\begin{equation}\label{equation} \sup\limits_{t \in \mathbb{R}} |\underline{x}^{'}(t)| \lesssim 1, \end{equation} 
\[ (u(x-\underline{x}(t),y,t) :t \in \mathbb{R})\quad \textmd{is precompact in} \quad H^1(\mathbb{R}^2 \times \mathbb{T}).
\]
\end{corollary}

\begin{theorem}\label{theorem}Assume that $u$ satisfies the conclusion of Corollary 7.5, then $u=0$. 
\end{theorem}
\noindent The proofs of the above three theorems are similar to those for the cubic case and we omit them. We also refer to [14, 19, 38]. \vspace{3mm}

\noindent \textbf{Acknowledgments.} I would like to express my deep thanks to my thesis advisor Professor Benjamin Dodson for many useful discussions, suggestions and comments. Also, I really appreciate Professor Zaher Hani, Professor Benoit Pausader, Qingtang Su and Professor Lifeng Zhao for insightful suggestions and valuable comments. Moreover, part of this work was progressed when the author attended the PDE conference (Waves, Spectral Theory and Applications Part 2) at UNC and I would like to thank the organizers for hosting.

\hfill \linebreak
\noindent \author{Zehua Zhao}

\noindent \address{Johns Hopkins University, Department of Mathematics, 3400 N. Charles Street, Baltimore, MD 21218, U.S.}

\noindent \email{zzhao25@jhu.edu}\\
\end{document}